\pdfoutput=1

\documentclass[12pt]{article}
\usepackage{graphicx}
\usepackage{amssymb, cite}
\usepackage{setspace} 
\usepackage[hmargin=1.15in,vmargin=1.15in]{geometry} 

\newcommand{\be}{\begin{equation}}
\newcommand{\ee}{\end{equation}}
\newcommand{\bqn}{\begin{eqnarray}}

\newcommand{\eqn}{\end{eqnarray}}

\newcommand{\bd}{\begin{description}}
\newcommand{\ed}{\end{description}}

\newtheorem{stat}{}[section]

\def\bs{\begin{stat}}
\def\es{\end{stat}}
\def\ben{\begin{enumerate}}
\def\een{\end{enumerate}}

\def\bp{\noindent{\bf Proof}  \ \ \ }

\newcommand{\ep}{\hfill $\square$}

\makeatletter
\@addtoreset{equation}{section}

\makeatother

\begin{document}

\begin{center}
{\large {\bf PACKING 3-VERTEX PATHS }}
\\[1ex]
{\large {\bf IN CUBIC 3-CONNECTED GRAPHS}}
\\[4ex]
{\large {\bf Alexander Kelmans}}
\\[2ex]
{\bf University of Puerto Rico, San Juan, Puerto Rico}
\\[0.5ex]
{\bf Rutgers University, New Brunswick, New Jersey}
\\[2ex]
\end{center}

\begin{abstract}
Let $v(G)$ and $\lambda (G)$ be the number 
of vertices and the maximum number of disjoint
3-vertex paths in $G$, respectively.
We discuss the following old
\\[0.5ex]
\indent
{\bf Problem.}
{\em Is the following claim true ?
\\
$(P)$ if $G$ is a 3-connected and cubic graph,
then $\lambda (G) =   \lfloor v(G)/3 \rfloor $.}
\\[0.5ex]
\indent
We show, in particular, that claim $(P)$  
is equivalent to some seemingly 
stronger claims (see {\bf \ref{3-con}}).
It follows that if claim $(P)$ is true, then
Reed's dominating graph conjecture is true for cubic 
3-connected graphs.
\\[1ex]
 
\end{abstract}

\section{Introduction}

\indent

We consider undirected graphs with no loops and 
no parallel edges. All notions and facts on graphs, that are  
used but not described here, can be found in \cite{BM,D,Wst}.
\\[1ex]
\indent
Given graphs $G$ and $H$, 
an $H$-{\em packing} of $G$ is a subgraph of $G$ 
whose components are isomorphic to $H$.
An $H$-{\em packing} $P$ of $G$ is called 
an $H$-{\em factor} if $V(P) = V(G)$. 
The $H$-{\em packing problem}, i.e. the problem of 
finding in $G$ an $H$-packing, having the maximum 
number of vertices, turns out to be $NP$-hard if $H$ is 
a connected graph with at least three vertices \cite{HK}.
Let $\Lambda $ denote a 3-vertex path.
In particular, the $\Lambda $-packing problem
is $NP$-hard.
Moreover, this problem  remains 
$NP$-hard even for cubic graphs \cite{K1}.

Although the $\Lambda $-packing problem is $NP$-hard, 
i.e. possibly intractable in general, this problem turns out 
to be tractable for some natural classes of graphs. 
It would be  also interesting to find polynomial algorithms 
that would provide a good approximation solution for 
the problem. Below (see {\bf \ref{km}}, {\bf \ref{eb(G)clfr}},  
and {\bf \ref{2conclfr}})
are some examples of such results.
In each case the corresponding packing problem is 
polynomially solvable.
\\[1ex]
\indent
Let $v(G)$ and $\lambda (G)$ denote the number 
of vertices and the maximum number of disjoint
3-vertex paths in $G$, respectively.
Obviously $\lambda (G) \le \lfloor v(G)/3 \rfloor $.

A graph is called {\em claw-free} if it contains no induced subgraph isomorphic
to $K_{1,3}$ (which is called a {\em claw}). 
A block of a connected graph is called an 
{\em end-block} if it has at most  one vertex in common 
with any other block of the graph. 
Let $eb(G)$ denote the number of end-blocks of $G$.
\bs {\em \cite{KKN}}
\label{eb(G)clfr}
Suppose that $G$ is a connected claw-free graph and  
$eb(G) \ge 2$. 
Then  $\lambda(G) \ge \lfloor (v(G) - eb(G) + 2)/3 \rfloor$,
and this lower bound is sharp.
\es

\bs {\em \cite{KKN}}
\label{2conclfr} 
Suppose that $G$ is a connected and claw-free graph 
having at most two end-blocks 
$($in particular, a 2-connected and claw-free graph$)$.
Then $\lambda(G) = \lfloor v(G)/3 \rfloor$.
\es 

Obviously the claim in {\bf \ref{2conclfr}} on claw-free graphs with exactly two end-blocks follows from  
{\bf \ref{eb(G)clfr}}.
\\[1ex]
\indent
In \cite{K,KM} we answered the following natural question:
\\[1ex]
\indent
{\em How many disjoint 3-vertex paths must a cubic graph have?}
\bs 
\label{km} If $G$ is a cubic graph then 
$\lambda (G) \ge \lceil v(G)/4  \rceil$.
Moreover, there is a polynomial time algorithm for 
finding a $\Lambda $-packing having at least  
$\lceil v(G)/4  \rceil$ components.
\es

Obviously if every component of $G$ is $K_4$, then
$\lambda (G) = v(G)/4$. 
Therefore the bound in {\bf \ref{km}} is sharp.
\\[1ex]
\indent
Let ${\cal G}^3_2 $ denote the set of graphs with each vertex of degree at least $2$ and at most $3$.

In \cite{K} we  answered (among other results) the following  question:
\\[1ex]
\indent
{\em How many disjoint 3-vertex paths must an $n$-vertex  graph from ${\cal G}^3_2$ have?}
\bs 
\label{2,3-graphs} Suppose that $G \in {\cal G}^3_2$ and  
$G$ has no 5-vertex components.
Then $\lambda (G) \ge \lceil v(G)/4  \rceil $.
\es

Obviously {\bf \ref{km}}  follows from {\bf \ref{2,3-graphs}} because if $G$ is a cubic graph, then $G \in {\cal G}^3_2$ 
and $G$ has no 5-vertex components. 
\\[1ex]
\indent
In \cite{K} we  also gave a construction  that allowed  to prove the following:
\bs 
\label{extrgraphs1}
There are infinitely many connected graphs for  which the bound  in {\bf \ref{2,3-graphs}} is attained.
Moreover, there are infinitely many subdivisions of
cubic 3-connected graphs for which the bound  in 
{\bf \ref{2,3-graphs}} is attained.
\es

The next interesting question is:
\\[1ex]
\indent
{\em How many disjoint 3-vertex paths must a cubic connected graph have?}
\\[1ex]
\indent
In \cite{K2} we proved the following.
\bs 
\label{cubic-connected} 
Let ${\cal C}(n)$ denote the set of connected cubic graphs with $n$ vertices and 
\\ 
$\lambda _n = \min \{\lambda (G)/v(G): G \in {\cal C}(n)\}$.
Then for some $c > 0$,
\[\frac{3}{11}(1 - \frac{c}{n}) \le \lambda _n \le 
\frac{3}{11}(1 - \frac{1}{n^2}).\]
\es

The next natural question is:
\bs {\bf Problem}
\label{2con-must}
How many disjoint 3-vertex paths must a cubic 2-connected graph have?
\es

This question is still open  (namely, the sharp lower bound on the number  of disjoint 3-vertex paths in a cubic 2-connected $n$-vertex graph is unknown).
\\[1ex]
\indent
On the other hand, it is also natural to consider the following 
\bs {\bf Problem.}
\label{Pr2con}
Are there 2-connected cubic graphs $G$ such that
$\lambda (G) <   \lfloor v(G)/3 \rfloor $ ?
\es

In \cite{K2con-cbp}
we gave a construction that provided infinitely many   
2-connected, cubic, bipartite, and planar  graphs such that  
$\lambda (G) <   \lfloor v(G)/3 \rfloor $.
\\[1ex]
\indent
The main goal of  this paper (see also \cite{K3con-cub}) is to  discuss the following old open problem which is similar to Problem {\bf \ref{Pr2con}}.
\bs {\bf Problem.}
\label{Pr3con} Is the following claim true ?
\\[1ex]
$(P)$ if $G$ is a 3-connected and cubic graph,
then $\lambda (G) =   \lfloor v(G)/3 \rfloor $.
\es

We show, in particular, that claim $(P)$ in 
{\bf \ref{Pr3con}} is equivalent to some seemingly 
stronger claims (see {\bf \ref{3-con}}).

In Section \ref{constructions} we give some notation, constructions, and simple observations.

In Section \ref{3connected} we formulate and prove our main theorem {\bf \ref{3-con}} concerning various claims  that are equivalent to claim $(P)$ in {\bf \ref{Pr3con}}.
We actually give different proofs of  {\bf \ref{3-con}}.
Thus if there is  a counterexample $C$ to one of the above claims, 
then the different proofs below provide different constructions 
of counterexamples to the other claims in {\bf \ref{3-con}}.
Moreover, different proofs provide better understanding of relations between various $\Lambda $-packing properties considered in {\bf \ref{3-con}}.

In \cite{R} B. Reed conjectured that
if $G$ is a connected cubic graph, then 
$\gamma (G) \le \lceil v(G)/3 \rceil $, where $\gamma (G)$ is the dominating number of $G$ (i.e. the size of a minimum vertex subset $X$ in $G$ such that every vertex in $G - X$ is adjacent to a vertex in $X$). It turns out that Reed's conjecture is not true for connected and even for 2-connected cubic graphs \cite{Kcntrex,KS}. 
If claim $(P)$ in {\bf \ref{Pr3con}} is true, then from 
{\bf \ref{3-con}} it follows, in particular, that Reed's conjecture is  true for 3-connected cubic graphs.

In Section \ref{AlmostCubic} we describe some  results showing that certain claims in {\bf \ref{3-con}} are best possible.

In Section \ref{homomorphism} we give a  $\Lambda $-factor homomorphism theorem in cubic graphs.

\section{Notation, constructions, and simple observations}
\label{constructions}

\indent
  
\begin{figure}
  \centering
  \includegraphics{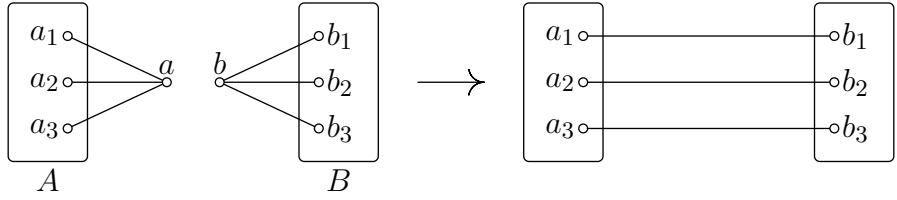}
  \caption{$Aa \sigma bB$}
  \label{fAasbB}
\end{figure}

We consider undirected graphs with no loops and 
no parallel edges unless stated otherwise. 
As usual, $V(G)$ and $E(G)$ denote the set of vertices and edges of $G$, respectively, and $v(G) = |V(G)$.
If $X$ is a vertex subset or a subgraph of $G$, then let
$D(X,G)$ or simply $D(X)$, denotes the set of edges in $G$, having exactly one end--vertex in $X$, and let
$d(X,G) = |D(X,G)|$.
If $x \in V(G)$, then $D(x,G)$ is the set of edges in $G$ incident to $x$,
$d(x,G) = |D(x,G)|$, $N(x,G) = N(x)$ is the set of vertices 
in $G$ adjacent to $x$, and 
$\Delta (G) = \max \{d(x,G): x \in V(G)\}$.
If $e = xy \in E(G)$, then let $End(e) = \{x,y\}$.
Let $Cmp (G)$ denote the set of components of $G$ and $cmp(G) = |Cmp(G)|$.
\\[1ex]
\indent
Let $A$ and $B$ be disjoint graphs, 
$a \in V(A)$,  $b \in V(B)$, and
$\sigma : N(a,A) \to N(b,B)$ be a bijection.
Let $Aa \sigma bB$ denote the graph
$(A - a) \cup (B - b) \cup \{x\sigma (x): x \in N(a,A)\}$.
We usually assume that 
$N(a,A) = \{a_1, a_2, a_3\}$, $N(b,B) = \{b_1, b_2, b_3\}$, and  $\sigma (a_i) = b_i$ for $i \in \{1,2,3\}$
(see Fig. \ref{fAasbB}). 
We also say that {\em $Aa \sigma bB$ is obtained from $B$ 
by replacing vertex $b$ by $(A - a)$ according to $\sigma $}.

Let $B$ be a cubic graph and $X \subseteq V(B)$. 
Let $A(v)$, where $v \in X$, be a graph, 
$a^v$ be  a vertex of degree three in $A(v)$, and 
$A^v = A(v) - a^v$.
By using the above operation, we can build a graph
$G = B\{(A(v), a^v): v \in X\}$ by replacing each vertex $v$
of $B$ in $X$ by $A^v $ assuming that all $A(v)$'s are disjoint. Let $D^v = D(A^v,G)$.
For each $u \in V(B) \setminus X$ let $A(u)$ be the graph having exactly two vertices $u$, $a^u$ and exactly three parallel edges connecting $u$ and $a^u$.
Then $G = B\{(A(v), a^v): v \in X\} = 
B\{(A(v), a^v): v \in V(B)\}$.
If, in particular, $X = V(G)$ and each $A(v)$ is a copy of $K_4$, then $G$ is obtained from $B$ by replacing each vertex by a triangle.

Let $E' = E(G) \setminus \cup \{E(A^v):  v \in V(B)\}$.
Obviously, there is a unique bijection
$\alpha : E(B) \to E'$ such that if $uv \in E(B)$, then 
$\alpha (uv)$ is an edge in $G$ having one end-vertex in
$A^u$ and the other in $A^v$.

Let $P$ be a $\Lambda $-packing in $G$.
For $uv \in E(B)$, $u \ne v$, we write  $u \neg ^p v$ or simply,
$u \neg  v$,
if $P$ has a 3-vertex path $L$ such that $\alpha (uv) \in E(L)$ and $|V(A^u) \cap V(L)| = 1$.
Let $P^v$ be the union of components of $P$ that meet
$D^v$ in $G$.
\\[1ex]
\indent
Obviously
\bs
\label{AasbB} 
Let $k$ be an integer and $k \le 3$.
If $A$ and $B$ above are $k$-connected, cubic, 
bipartite, and planar graphs, then $Aa\sigma bB$ 
is also a $k$-connected, cubic, bipartite, and 
planar graph, respectively.
\es

From {\bf \ref{AasbB}} we have: 
\bs
\label{A(B)} 
Let $k$ be an integer and $k \le 3$.
If $B$  and each $A_v$ is a $k$-connected, cubic, 
bipartite, and planar graphs, then
$B\{(A_v, a_v): v \in V(B)\}$ is also a $k$-connected, 
cubic, bipartite, and planar graph, respectively.
\es

\begin{figure}
  \centering
  \includegraphics{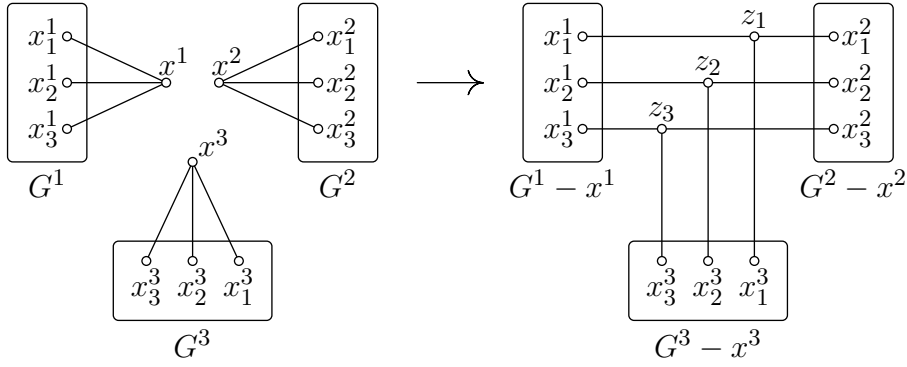}
  \caption{$Y(A^1,a^1; A^2,a^2; A^3,a^3)$}
  \label{fY}
\end{figure}

Let $A^1$, $A^2$, $A^3$ be three disjoint graphs,
$a^i \in V(A^i)$,  and $N(a^i,A^i) = \{a^i_1, a^i_2, a^i_3\}$, 
where $i \in \{1,2,3\}$.
Let  $F = Y(A^1,a^1; A^2,a^2; A^3,a^3)$  denote the graph 
obtained from $(A^1 - a^1) \cup (A^2-a^2) \cup (A^3 - a^3)$ 
by adding three new vertices $z_1$, $z_2$, $z_3$ 
and the set of nine new edges $\{z_ja^i_j:  i, j \in \{1,2,3\}$
(see Fig. \ref{fY}). 
In other words, if $B = K_{3,3}$ is the complete
$(X,Z)$-bipartite graph with 
$X = \{x_1,x_2,x_3\}$ and  $Z = \{z_1,z_2,z_3\}$, then $F$ is obtained from the $B$ by replacing each vertex $x_i$ in $X$ by $A^i - a^i$ so that
$D(A^i - a^i, F) = \{a^i_jz_j: j \in \{1,2,3\}$.
Let $D^i = D(A^i- a^i, F)$.
If $P$ is a $\Lambda $-packing of $F$, then let 
$P^i = P^i(F)$ be the union of components of $P$ meeting $D^i$ and $E^i = A^i(P) = E(P) \cap D^i$, $i \in \{1,2,3\}$.

If each $(A^i, a^i)$ is a copy of the same $(A, a)$, 
then we write $Y(A, a)$ instead of 
$Y(A^1, a^1; A^2,a^2; A^3,a^3)$.
\\[1.5ex]
\indent
From {\bf \ref{A(B)} } we have, in particular: 
\bs
\label{Y(G1,G2,G3)} 
Let $k$ be an integer and $k \le 3$.
If each $A^i$ above is a $k$-connected, cubic, and 
bipartite graph, then $Y(A^1,a^1; A^2,a^2; A^3,a^3)$ 
{\em (see Fig. \ref{fY})} is also a $k$-connected, cubic, and bipartite graph, 
respectively.
\es

We will use the following simple observation.
\bs
\label{3cut}
Let $A$ and $B$ be disjoint graphs, 
$a \in V(A)$,  $N(a,A) = \{a_1, a_2, a_3\}$,  
$b \in V(B)$,  $N(b,B) = \{b_1, b_2, b_3\}$, and 
$G = Aa \sigma bB$, where each $\sigma (a_i) = b_i$
{\em (see Fig. \ref{fAasbB})}. 
Let $P$ be a $\Lambda $-factor of $G$ 
$($and so $v(G) = 0 \bmod 3$$)$ and 
$P'$ be the $\Lambda $-packing of $G$ consisting of the components $($3-vertex paths$)$ of $P$ that meet 
$\{a_1b_1, a_2b_2, a_3b_3\}$.
\\[1ex]
$(a1)$
Suppose that  $v(A) = 0 \bmod 3$, and 
so $v(B) = 2 \bmod 3$.  Then one of the following holds
{\em (see Fig \ref{A1})}:
\begin{figure}
  \centering
  \includegraphics{lambda-packing-3.mps}
  \caption{}
  \label{A1}
\end{figure}
\\[0.5ex]
\indent
$(a1.1)$ $P'$ has exactly one component that
has two  vertices in $A - a$, that are adjacent 
$($and, accordingly, exactly one vertex in $B - b$$)$,
\\[0.5ex]
\indent
$(a1.2)$ $P'$ has exactly two components  and each component has exactly one vertex in $A - a$ 
$($and, accordingly, exactly two vertices in 
$B - b$, that are adjacent$)$,
\\[0.5ex]
\indent
$(a1.3)$ $P'$ has exactly three components $L_1$, $L_2$, $L_3$ and one of them, say $L_1$, has exactly one vertex in $A - a$ and each of the other two $L_2$, $L_3$, has exactly two vertices in $A - a$, that are adjacent
$($and, accordingly, $L_1$ has exactly two  vertices in $B - b$, that are adjacent, and each of the other two $L_2$, $L_3$, has exactly  one vertex in $B - b$$)$. 
\\[1ex]
$(a2)$
Suppose that  $v(A) = 1 \bmod 3$, and 
so $v(B) = 1 \bmod 3$.  Then one of the following holds
{\em (see Fig \ref{A2})}:
\begin{figure}
  \centering
  \includegraphics{lambda-packing-4.mps}
  \caption{}
  \label{A2}
\end{figure}
\\[0.5ex]
\indent
$(a2.1)$ $P' = \emptyset $,
\\[0.5ex]
\indent
$(a2.2)$ $P'$ has exactly two components, say $L_1$, $L_2$, and one of the them, say $L_1$, has exactly one vertex in 
$A - a$ and exactly two  vertices in $B - b$, that are adjacent, and the other component $L_2$ has exactly two vertices in 
$A - a$, that are adjacent, and exactly one vertex in $B - b$,
\\[0.5ex]
\indent
$(a2.3)$ $P'$ has exactly three components $L_1$, $L_2$, $L_3$ and either each $L_i$ has exactly one vertex in 
$A - a$ 
$($and, accordingly, has exactly two vertices in 
$B - b$, that are adjacent$)$ or each $L_i$ has exactly two vertices in $A - a$, that are adjacent 
$($and, accordingly, has exactly one vertex in $B - b$$)$.
\es

\section{$\Lambda $-packings in cubic 3-connected graphs}
\label{3connected}

\indent

The main goal of this section is to prove the following theorem
showing that claim $(P)$ in 
{\bf \ref{Pr3con}} is equivalent to various seemingly 
stronger claims.
\bs
\label{3-con}
The following are equivalent for cubic 3-connected graphs $G$:
\\[1ex]
${\bf (z1)}$ 
$v(G) = 0 \bmod 6$ $\Rightarrow$ $G$ has 
a $\Lambda $-factor,
\\[1ex] 
${\bf (z2)}$
$v(G) = 0 \bmod 6$ $\Rightarrow$ for every 
$e \in E(G)$ there is a $\Lambda $-factor of $G$ 
avoiding $e$ $($i.e. $G - e$ has a $\Lambda $-factor$)$,
\\[1ex]
${\bf (z3)}$
$v(G) = 0 \bmod 6$ $\Rightarrow$ for every 
$e \in E(G)$ there is a $\Lambda $-factor of $G$ 
containing $e$,
\\[1ex] 
${\bf (z4)}$
$v(G) = 0 \bmod 6$ $\Rightarrow$ for every 
$x \in V(G)$ there is at least one 3-vertex path $L$ 
such that $L$ is centered at $x$ and $G - L$ has 
a $\Lambda $-factor, 
\\[1ex]
${\bf (z5)}$
$v(G) = 0 \bmod 6$ $\Rightarrow$ for every 
$x \in V(G)$ there is at least two 3-vertex paths $L$ 
such that $L$ is centered at $x$ and $G - L$ has 
a $\Lambda $-factor, 
\\[1ex]
${\bf (z6)}$
$v(G) = 0 \bmod 6$ $\Rightarrow$ for every 
$xy \in E(G)$ there are edges $xx', yy' \in E(G)$ such that
$G - xyy'$ and $G - x'xy$ have $\Lambda $--factors,
\\[1ex]
${\bf (z7)}$
$v(G) = 0 \bmod 6$ $\Rightarrow$ $G - X$ has 
a $\Lambda $-factor for every $X \subseteq E(G)$ such that  $|X| = 2$,
\\[1ex]  
${\bf (z8)}$
$v(G) = 0 \bmod 6$ $\Rightarrow$ 
$G - L$ has a $\Lambda $-factor for every
3-vertex path $L$ in $G$,
\\[1ex]  
${\bf (z9)}$
$v(G) = 0 \bmod 6$ $\Rightarrow$ for every 3-edge cut $K$ of $G$ and $S \subset K$, $|S| = 2$, there is a 
$\Lambda $-factor $P$ of $G$ such that $E(P) \cap K = S$,
\\[1ex]
${\bf (t1)}$
$v(G) = 2 \bmod 6$ $\Rightarrow$ for every  
$x \in V(G)$ there is $xy \in E(G)$ such that $G - \{x,y\}$ 
has a $\Lambda $-factor,
\\[1ex]
${\bf (t2)}$
$v(G) = 2 \bmod 6$ $\Rightarrow$ $G - \{x,y\}$ 
has a $\Lambda $-factor for every $xy \in E(G)$,
\\[1ex]
${\bf (t3)}$
$v(G) = 2 \bmod 6$ $\Rightarrow$ for every 
$x \in V(G)$ there is a 5-vertex path $W$ such that $x$ is 
the center vertex of $W$ and $G - W$ has 
a $\Lambda $-factor 
{\em (see also {\bf \ref{3cut}} $(a1.2)$ and Fig \ref{A1})},
\\[1ex]
${\bf (t4)}$
$v(G) = 2 \bmod 6$ $\Rightarrow$ for every 
$x \in V(G)$ and $xy \in E(G)$ there is a 5-vertex path $W$
such that $x$ is the center vertex of $W$,
$xy \not \in E(W)$,  and $G - W$ has a $\Lambda $-factor
{\em (see also {\bf \ref{3cut}} $(a1.2)$ and Fig \ref{A1})},
\\[1ex]
${\bf (f1)}$
$v(G) = 4 \bmod 6$ $\Rightarrow$ $G - x$ 
has a $\Lambda $-factor for every $x \in V(G)$,
\\[1ex]
${\bf (f2)}$
$v(G) = 4 \bmod 6$ $\Rightarrow$ $G - \{x, e\}$ 
has a $\Lambda $-factor for every $x \in V(G)$ and 
$e \in E(G)$,
\\[1ex]
${\bf (f3)}$
$v(G) = 4 \bmod 6$ $\Rightarrow$ for every 
$x \in V(G)$ there is a 4-vertex path $Z$ such that
$x$ is an inner vertex of $Z$ and $G - Z$ has 
a $\Lambda $-factor
{\em (see also {\bf \ref{3cut}} $(a2.2)$ and Fig \ref{A2})},
\\[1ex]
${\bf (f4)}$
$v(G) = 4 \bmod 6$ $\Rightarrow$ for every 
$x \in V(G)$ there is $xy \in E(G)$ and a 4-vertex path $Z$ such that $x$ is an inner vertex of $Z$, $xy \not \in E(Z)$, 
and $G - Z$ has a $\Lambda $-factor 
{\em (see also {\bf \ref{3cut}} $(a2.2)$ and Fig \ref{A2})},
\\[1ex]
${\bf (f5)}$
$v(G) = 4 \bmod 6$ $\Rightarrow$ for every 
$xy \in E(G)$ there exists a 4-vertex path $Z$ such that $xy$ is the middle edge of $Z$ 
and $G - Z$ has a $\Lambda $-factor 
{\em (see also {\bf \ref{3cut}} $(a2.2)$ and Fig \ref{A2})},
\\[1ex]
${\bf (f6)}$
$v(G) = 4 \bmod 6$ $\Rightarrow$ for every $z \in V(G)$ and every 3-vertex path
$xyz$ there exists a 4-vertex path $Z$ such that $xyz \subset Z$, $z$ is an end-vertex of $Z$, 
and $G - Z$ has a $\Lambda $-factor 
{\em (see also {\bf \ref{3cut}} $(a2.2)$ and Fig \ref{A2})}.   
\es

Theorem {\bf \ref{3-con}} follows from 
{\bf \ref{z1ottoz2}} -- {\bf \ref{z9ottoz1}} below.

In \cite{Kclfree} we have shown that claims $(z1)$ - $(z5)$,
$(t1)$, $(t2)$, $(f1)$, and $(f2)$ are true for cubic, 
3-connected, and claw-free graphs.
\\[1ex]
\indent
The remarks below show that if claims $(z7)$, $(z8)$, $(t2)$,
$(f1)$, $(f2)$ in {\bf \ref{3-con}} 
are true, then they are best possible in some sense.
\\[0.5ex]
\indent
$(r1)$ Obviously claim $(z7)$ is not true if 
condition ``$|X| = 2$'' is replaced by condition ``$|X| = 3$''.
Namely, if 
$G$ is a cubic 3-connected graph,
$v(G) = 0 \bmod 6$, 
$X$ is a 3-edge cut in $G$, and the two components of 
$G - X$ have different number of vertices $\bmod 6$, 
then clearly $G - X$ has no  $\Lambda $-factor.
Also in Section \ref{AlmostCubic} (see {\bf \ref{CalY}}) we describe an infinite set of cubic 3-connected graphs $G$  having a triangle $T$ such that $G - E(T)$ has no  $\Lambda $-factor.
\\[0.5ex]
\indent
$(r2)$ There exist infinitely many triples $(G, L, e)$ such that
$G$ is a cubic, 3-connected, bipartite, and planar graph,
$v(G) = 0 \bmod 6$, $L$ is a 3-vertex path in $G$,
$e \in E(G - L)$, and $(G - e) - L \}$ has no 
$\Lambda $-factor, and so claim $(z8)$ is tight.
Moreover, there are infinitely many triples 
$(G, L, L')$ such that
$G$ is a cubic  3-connected graph, $L$ and $L'$ are disjoint 3-vertex paths in $G$, and $G - (L \cup L')$ has 
no $\Lambda $-factor. If $G$ has a triangle, then it is easy to find such $L$ and $L'$. Indeed,
let  $v \in V(G)$,
$N(v,G) = \{x,y,z\}$, and $yz \in E(G)$, and so $vyz$ is a triangle. Since $G$ is 3-connected, $x$ is not adjacent to $\{y,z\}$. Let $L$ and $L'$ be 3-vertex paths in $G - v$
containing $x$ and $yz$, respectively.
Then $v$ is an isolated vertex in $G - (L \cup L')$, and so
$G - (L \cup L')$ has  no $\Lambda $-factor.
Similar idea can be used to find such $L$ and $L'$ if $G$ has a 4-cycle. In Section \ref{AlmostCubic} (see {\bf \ref{G,L,L'}}) we describe  a sequence of 
infinitely many triples $(G, L, L')$ with the above property,
where $G$ is a cubic cyclically 6-connected graph,  and so
$G$ has no triangles, no 4-cycles,  and 5-cycles.
Thus  claim $(z8)$ is tight in this sense as well.
\\[0.5ex]
\indent
$(r3)$ 
There exist infinitely many triples $(G, xy, e)$ such that
$G$ is a cubic, 3-connected, bipartite, and planar 
graph, $v(G) = 2 \bmod 6$, $xy \in E(G)$, 
$e \in E(G - \{x,y\})$,
and $G - \{x, y ,e \}$ has no $\Lambda $-factor, and 
so claim $(t2)$ is tight.
\\[0.5ex]
\indent
$(r4)$ There exist infinitely many triples $(G, x, y)$ such that
$G$ is a cubic 3-connected graph,
$v(G) = 2 \bmod 6$, $\{x,y\} \subset V(G)$, $x \ne y$, 
$xy \not \in E(G)$, and $G -  \{x,y\}$ has no 
$\Lambda $-factor, and so 
claim $(t2)$ is  not true if vertices $x$ and $y$ 
are not adjacent.
\\[0.5ex ]
\indent
$(r5)$ 
There exist infinitely many $(G,a,b,x)$ such that
$G$ is a cubic, 3-connected graph with no 3-cycles and no 4-cycles, $v(G) = 4 \bmod 6$, $x \in V(G)$,
$a$ and $b$ are non-adjacent edges in $G - x$,
and $G - \{x, a,b \}$ has no $\Lambda $-factor, and 
so claim $(f2)$ is tight.
\\[0.5ex] 
\indent
$(r6)$ 
There exist  infinitely many triples $(G,L,x)$ such that $G$ 
is a cubic, 3-connected graph 
with no 3-cycles and no 4-cycles, $v(G) = 4 \bmod 6$,
$x \in V(G)$, $L$ is a 3-vertex path in $G - x$, and 
$G - \{x, L \}$ has no $\Lambda $-factor (see claim $(f1)$).
\\[1ex]
\indent 
 We need the following two results obtained before.

\bs {\em \cite{K2con-cbp}}
\label{Y,cmp(Pi)<3} 
Let $G = Y(A^1,a^1;  A^2,a^2; A^3,a^3)$ 
{\em (see Fig. \ref{fY})} and 
$P$  be a $\Lambda $-factor of $G$.
Suppose that each $A^i$ is a cubic graph and 
$v(A^i) = 0 \bmod 6$.
Then $cmp (P^i) \in \{1, 2\}$ for every $i \in \{1,2,3\}$.
\es

\bp Let $i \in \{1,2,3\}$.
Since $D^i$ is a matching and $P^i$ consists of the components of $P$ meeting $D^i$, clearly $cmp (P^i) \le 3$.
Since $v(A^i) = -1 \bmod 6$, we have $cmp(P^i) \ge 1$.
It remains to show that $cmp (P^i) \le 2$.
Suppose, on the contrary, that $cmp(P^1) = 3$.
Since $P$ is a $\Lambda $-factor of $G$ and 
 $v(A^1 - a^1) = -1\bmod 6$, clearly
 $v(P^1) \cap V(A^1 - a^1) = 5$ and we can assume 
 (because of symmetry) that 
$P_1$ consists of three components  $a^1_3z_3a^2_3$,  
$z_1a^1_1 y^1$, and $z_2a^1_2 u^1$ for some 
$y^1, u^1 \in V(A^1)$ 
Then $cmp(P^3) = 0 $, a contradiction.
\ep
\bs {\em \cite{K2con-cbp}}
\label{Y,e-} 
Let $A$ be a graph, $e  = aa_1\in E(A)$, and 
$G = Y(A,a)$ {\em (see Fig. \ref{fY})}.
Suppose that
\\[0.5ex]
$(h1)$ $A$ is  cubic,
\\[0.5ex]
$(h2)$ $v(A) = 0\bmod 6$, and 
\\[0.5ex]
$(h3)$ $a$ has no $\Lambda $-factor containing 
$e = aa_1$.
\\[0.5ex]
\indent
Then $v(G) = 0 \bmod 6$ and $G$ has 
no $\Lambda $-factor.
\es

\bp (uses {\bf \ref{Y,cmp(Pi)<3}}).
Suppose, on the contrary, that $G$ has  a $\Lambda $-factor $P$. By definition of $G = Y(A,a)$, each $A^i$ is a copy of $A$ and edge $e^i = a^ia^i_1$ in $A^i$ is a copy of edge 
$e = aa_1$ in $A$ 
By {\bf \ref{Y,cmp(Pi)<3}}, $cmp(P^i) \in \{1, 2\}$
for every $i \in \{1,2,3\}$.
Since $P$ is a $\Lambda $-factor of $G$ and 
$v(A^i - x^i) = -1\bmod 6$, clearly
$E(P) \cap D^i $ is an edge subset of a $\Lambda $--factor 
of $A^i$ for every $i \in \{1,2,3\}$
(we assume that edge $z_ja^i_j$ in $G$ is edge $a^ia^i_j$ in $A^i$).
Since $a^1a^i_1$ belongs to no  $\Lambda $-factor of 
$A^i$ for every $i \in \{1,2,3\}$, clearly
$E(P) \cap \{z_1a^1_1,z_1a^2_1,z_1a^3_1\}  = \emptyset $.
Therefore $z_1 \not \in V(P)$, and so $P$ is not 
a $\Lambda $-factor of $G$, a contradiction.
\ep 
\bs 
\label{z1ottoz2}
$(z1)$ $\Leftrightarrow$ $(z2)$.
\es
 
\bp (uses {\bf \ref{Y(G1,G2,G3)}} and {\bf \ref{Y,cmp(Pi)<3}}). 
Obviously $(z1)$ $\Leftarrow$ $(z2)$. 
We prove $(z1)$ $\Rightarrow$ $(z2)$.

Suppose, on the contrary, that $(z1)$ is true but $(z2)$ 
is not true, i.e. there is a cubic 3-connected graph $A$ 
and $aa_1 \in E(G)$  such that $v(A) = 0 \bmod 6$ and
every $\Lambda $-factor of $G$ contains $aa_1$.
Let $G = Y(A^1,a^1; A^2, a_2; A^3, a^3)$, where 
each $(A^i, a^i)$ above is a copy of $(A, a)$ and 
edge $a^ia^i_1$ in $A^i$ is a copy of edge $aa_1$ 
in $A$ 
(see Fig. \ref{fY}).
Since $A$ is cubic and 3-connected, 
by {\bf \ref{Y(G1,G2,G3)}},  $G$ is also cubic and 
3-connected.
Obviously $v(G) = 0  \bmod 6$.
By $(z1)$, $G$ has a $\Lambda $-factor $P$.
Since each $v(A^i) = 0\bmod 6$, by {\bf \ref{Y,cmp(Pi)<3}},
$cmp(P^i)  \in \{1, 2\} $ for every $i \in \{1,2,3\}$.
Since $P$ is a $\Lambda $-factor of $G$ and 
$v(A^i - a^i) = -1\bmod 6$, clearly
$E(P) \cap D^i $ is an edge subset of 
a $\Lambda $-factor of $A^i$ for every $i \in \{1,2,3\}$
(we assume that edge $z_ja^i_j$ in $G$ is edge $a^ia^i_j$ 
 in $A^i$).
Since $a^1a^i_1$ belongs to every  $\Lambda $-factor of 
$A^i$ for every $i \in \{1,2,3\}$, clearly
$ \{z_1a^1_1, z_1a^2_1,z_1a^3_1\} \subseteq E(P)$.
Therefore vertex $z_1$ has degree three in $P$, and so 
$P$ is not a $\Lambda $-factor of $G$, a contradiction.
\ep
\bs 
\label{z1ottoz3}
$(z1)$ $\Leftrightarrow$ $(z3)$.
\es

\bp Claim $(z1) \Leftarrow (z3)$ is obvious. 
Claim $(z1)\Rightarrow (z3)$ follows from {\bf \ref{Y,e-}}.
\ep
\\

Let $B$ be a cubic graph. Given $v \in V(B)$,
let $A(v)$ be a cubic graph, 
$a^v \in V(A(v))$, and $A^v = A(v) - a^v$.
We assume that all $A(v)$'s are disjoint.
Let $G$ be a graph obtained from $B$ by replacing each vertex $v$ in $B$ by $A^v$.
Let $D^v = D(A^v,G)$ and $E' = E(G) \setminus \cup \{E(A^v):  v \in V(B)\}$.
Obviously, there is a bijection
$\alpha : E(B) \to E'$ such that if $uv \in E(B)$, then 
$\alpha (uv)$ is an edge in $G$ having one end-vertex in
$A^u$ and the other in $A^v$.

Let $P$ be a $\Lambda $-packing in $G$.
For $uv \in E(B)$, $u \ne v$, we write  $u \neg v$
if $P$ has a component $L$, such that $\alpha (uv) \in E(L)$ and $|V(A^u) \cap V(L)| = 1$.
Let $P^v$ be the union of components of $P$ that meet
$D^v$ in $G$.
\bs 
\label{z1ottoz4}
$(z1)$ $\Leftrightarrow$ $(z4)$.
\es

\bp (uses {\bf \ref{A(B)}} and {\bf \ref{3cut} }$(a1)$). 
Obviously $(z1) \Leftarrow (z4)$.
We prove $(z1) \Rightarrow (z4)$.

\begin{figure}
  \centering
  \includegraphics{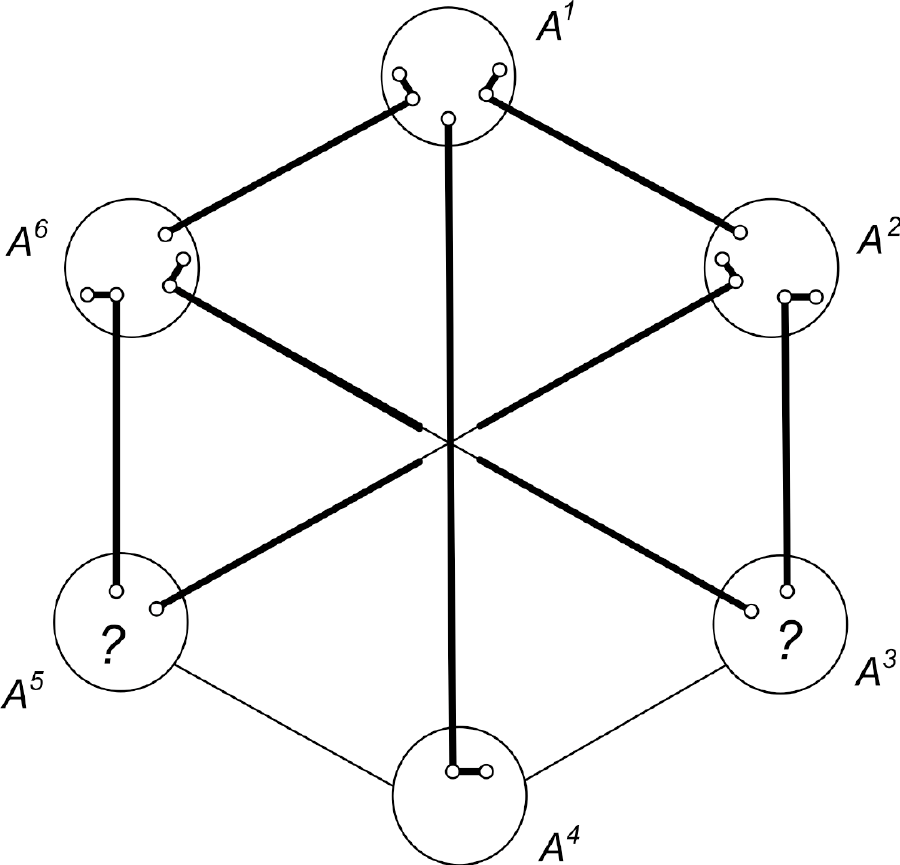}
  \caption{$(z1)$ $\Rightarrow$ $(z4)$}
  \label{fPrz1toz4}
\end{figure}

Suppose that $(z1)$ is true but  $(z4)$ is not true, 
i.e.  there is a cubic 3-connected graph $A$
and $a\in V(A)$  such that $v(A) = 0 \bmod 6$ and
$a$ has degree one in every $\Lambda $-factor of $A$.
Let $G$ be the graph obtained from $B = K_{3,3}$ by replacing each vertex $v$ by a copy $A^v$ of $A - a$ 
(see Fig. \ref{fPrz1toz4}). 
Obviously $v(G) = 0 \bmod 6$ and by {\bf \ref{A(B)}},
$G$ is a cubic, 3-connected graph.
By $(z1)$, $G$ has a $\Lambda $-factor $P$.
If $uv \in E(P)$, then let $L(uv)$ denote the component 
of $P$ containing $uv$. Let $V(B) = \{1, \ldots , 6\}$.

Since vertex $a$ has degree one in every $\Lambda $-factor of $A$ and each $(A(v), a^v)$ is a copy of $(A,a)$, by
{\bf \ref{3cut} }$(a1)$, we have: $cmp(P^v) \in \{1, 3\}$.
\\[1ex]
${\bf (p1)}$
Suppose that there is $v \in V(B)$ such that 
$cmp(P^v) = 3$. By symmetry of $B$, we can assume that
$v = 1$ and,  by {\bf \ref{3cut}} $(a1.3)$,
$1 \neg 4$, $6 \neg 1$, and $2 \neg 1$
(see Fig. \ref{fPrz1toz4}). 
Let $x \in \{2,6\}$. Since $cmp(P^x) \in \{1,3\}$ and
$|V(L(1x) \cap V(A^x)| = 1$, clearly $cmp(P^x) = 3$ and,  
by {\bf \ref{3cut}} $(a1.3)$, $5 \neg x$, $3 \neg x$. 
Then by {\bf \ref{3cut} }$(a1.2)$, 
$cmp(P^s) = 2$ for $s \in \{3,5\}$, a contradiction.
\\[1ex]
${\bf (p2)}$ Now suppose that $cmp(P^v) = 1$ for every 
$v \in V(B)$. By symmetry, we can assume
$1 \neg 2$. Then $P^1$ contradicts {\bf \ref{3cut}} $(a1.1)$.
\ep
\bs
\label{z1ottot1}
$(z1)\Leftrightarrow (t1)$.
\es

\bp (uses {\bf \ref{Y(G1,G2,G3)}})
Let $G = Y(A^1,a^1; A^2,a_2; A^3,a^3)$
(see Fig. \ref{fY}).
By {\bf \ref{Y(G1,G2,G3)}}, if each $A^i$ is cubic 
and 3-connected then $G$ is also cubic and 
3-connected.
\\[1ex]
${\bf (p1)}$
We first prove $(z1) \Rightarrow (t1)$.
Suppose, on the contrary, that $(z1)$ is true but
$(t1)$ is not true, i.e. there is a cubic 3-connected graph 
$A$ and $a \in V(A)$ such that $v(A) = 2 \bmod 6$ and
 $A - \{a,y\}$ has no $\Lambda $-factor for every vertex 
 $y$ in $A$ adjacent to $a$. 
Let each $(A^i, a^i)$ above be a copy of $(A, a)$, and so
$A^i - \{a^i,a^i_j\}$ has no $\Lambda $-factor for every
$i,j \in \{1,2,3\}$. 
Obviously $v(G) = 0 \bmod 6$. By $(z1)$, $G$ has 
a $\Lambda $-factor $P$.
Then it is easy to see that  since each 
$v(A^i - a^i) = 1 \bmod 6$, 
there are $r,s, j \in \{1,2,3\}$ such that $r \ne s$ and 
$P^s = x^s_jz_jx^r_j$. 
Since $P$ is a $\Lambda $-factor of $G$, clearly 
$P \cap (A^s - \{z_j,a^s_j\})$ is a $\Lambda $-factor 
of $A^s - \{z_j,a^s\} = A^s - \{a^s,x^s_j\}$, a contradiction. 
\\[1ex]
${\bf (p2)}$
Now we prove $(z1) \Leftarrow (t1)$.
Suppose, on the contrary, that $(t1)$ is true but
$(z1)$ is not true, i.e. there is a cubic 3-connected
graph $A$ such that $v(A) = 0 \bmod 6$ and 
$A$ has no $\Lambda $-factor. 

Let $(A^i,a^i)$ above be a copy of $(A,a)$ for 
$i \in \{1,2\}$, where $a \in V(A)$, and $(A^3,a^3)$ be 
a copy of $(H,h)$  for some cubic 3-connected graph 
$H$ and $h \in V(H)$, where $v(H) = 2 \bmod 6$. 
Obviously $v(G) = 2 \bmod 6$.
Suppose that $P$ is a $\Lambda $-factor of 
$G - \{z_1a^3_1\}$. Then $cmp(P^1) \le 2$. 
Since $v(A^1 - a^1) = -1 \bmod 6$, we have $cmp(P^1)\ge 1$.
Now since $P$ is a $\Lambda $--factor of 
$G- \{z_1a^3_1\}$ and $v(A^1- a^1) = -1\bmod 6$, 
clearly $E(P) \cap D^1 $ is an edge subset of 
a component of a $\Lambda $-factor of $A^1$
(we assume that edge $z_ja^1_j$ in $G$ is edge $a^1a^1_j$ in $A^1$).
Therefore $A$ has a $\Lambda $-factor, a contradiction.
\ep
\bs 
\label{z1ottof1}
$(z1)$ $\Leftrightarrow$ $(f1)$.
\es

\bp (uses {\bf \ref{Y(G1,G2,G3)}}) 
Let $G = Y(A^1,a^1; A^2,a_2; A^3,a^3)$
(see Fig. \ref{fY}).
By {\bf \ref{Y(G1,G2,G3)}}, if each $A^i$ is 
cubic and 3-connected then $G$ is also cubic and 
3-connected.
\\[1ex]
${\bf (p1)}$
We first prove $(z1) \Rightarrow (f1)$.
Suppose, on the contrary, that $(z1)$ is true but $(f1)$ 
is not true, i.e. there is a cubic 3-connected graph $A$  
such that $v(A) = 4 \bmod 6$ but $A - a$ has no 
$\Lambda $-factor for some $a \in V(A)$.
Let each $(A^i,a^i)$ above be a copy of $(A, a)$.
Obviously $v(G) = 0 \bmod 6$.
By $(z1)$, $G$ has a  $\Lambda $-factor $P$.
Let $E^i = E^i(P) =D^i \cap E(P)$.
 
Since $P$ is a  $\Lambda $-factor of $G$ and 
$v(A^i - a^i) = 3 \bmod 6$, clearly $|E^i| \in \{0, 2, 3\}$.
Since $A^i - a^i$ has no $\Lambda $-factor, 
$|E^i| \in \{2, 3\}$ for $i \in \{1,2,3\}$. 
Then each $d(z_j, P) \ge 2$. 
Since $P$ is a $\Lambda $-factor of $G$, 
clearly each $d(z_j, P) \le 2$.
Therefore each $d(z_j, P) = 2$. But then
$|E^i| = 1$ for some $i \in \{1,2,3\}$, a contradiction.
\\[1ex]
${\bf (p2)}$
Now we prove $(z1) \Leftarrow (f1)$.
Suppose, on the contrary, that $(f1)$ is true but $(z1)$ is 
not true, i.e. there is a cubic 3-connected graph $A$  
such that $v(A) = 0 \bmod 6$ and $A$ has no 
$\Lambda $-factor.
Let $(A^i,a^i)$ above be a copy of $(A, a)$ for
$i \in \{1,2\}$ and some $a \in V(A)$ and $(A^3,a^3)$ 
is a copy of $(H,h)$ for some cubic 3-connected graph 
$H$ and $h \in V(H)$, where $v(H) = 4 \bmod 6$ 
(see Fig. \ref{fY}).
Then $v(G) = 4 \bmod 6$. Let $x \in V(H - h)$.
Suppose that $G - x $  has a $\Lambda $-factor $P$.
Since $A$ has no $\Lambda $-factor, clearly 
$|E^1(P)| = |E^2(P)| = 3$. Then $P$ is not a 
$\Lambda $-factor of $G - x$, and so $(f1)$ is not true, 
a contradiction.
\ep
\bs 
\label{z4ottot2}
$(z4)$ $\Leftrightarrow$ $(t2)$.
\es

\bp (uses {\bf \ref{Y(G1,G2,G3)}}, {\bf \ref{z1ottoz4}}, and 
{\bf \ref{z1ottot1}}).
Let $G = Y(A^1,a^1; A^2,a_2; A^3,a^3)$
(see Fig. \ref{fY}).
By {\bf \ref{Y(G1,G2,G3)}}, if each $A^i$ is 
cubic and 3-connected, then $G$ is also cubic and 
3-connected.
\\[1ex]
${\bf (p1)}$
We first prove $(z4) \Leftarrow (t2)$.
Obviously $(t2)  \Rightarrow (t1)$.
By  {\bf \ref{z1ottoz4}},
$(z1)$ $\Leftrightarrow$ $(z4)$.
By  {\bf \ref{z1ottot1}}, $(z1)\Leftrightarrow (t1)$.
The result follows.
\\[1ex]
${\bf (p2)}$ Now we prove
$(z4)\Rightarrow (t2)$.
Suppose, on the contrary, that $(z4)$ is true but $(t2)$ 
is not true. Then there is a cubic 3-connected graph $A$ 
and $aa_1 \in E(G)$ such that $v(A) = 2 \bmod 6$ and 
$A - \{a,a_1\}$ has no $\Lambda $-factor.
Let each $(A^i, a^i)$ above be a copy of $(A,a)$ and 
edge $a^ia^i_1$ in $A^i$ be a copy of edge $aa_1$ in $G$.
Obviously $v(G) = 0 \bmod 6$.
Let $L^i = a^j_1z_1a^k_1$, where $\{i, j, k\} = \{1,2,3\}$.
By $(z4)$, $G$ has a $\Lambda $-factor $P$ containing 
$L^i$ for some $i \in \{1,2,3\}$, say for $i = 3$. 
If $s \in \{1,2\}$, then  $cmp (P^s) = 3$ because
$A_s - \{a^s, a^s_1\}$ has no $\Lambda $-factor and 
$v(A^s - a^s) = 1 \bmod 6$.
Also $cmp(P^3) \ge 1$ because
 $v(A^3 - a^3) = 1 \bmod 6$.
Then $P_1 \cup P^2 \cup P^3$ has at least four components each meeting $\{z_1,z_2,z_3\}$, 
a contradiction.
\ep 
\bs 
\label{z2toz5}
$(z2)$ $\Rightarrow$ $(z5)$.
\es

\bp (uses  {\bf \ref{AasbB}} and {\bf \ref{3cut}} $(a1)$).
Suppose, on the contrary, that $(z2)$ is true but $(z5)$ 
is not true. Then there is a cubic 3-connected graph $A$
and $a \in V(A)$ such that at most one 3-vertex path, centered at $a$ and belonging to a  $\Lambda $-factor of $A$. It is sufficient to prove our claim in case when
$A$ has exactly one 3-vertex path, say $L = a_1 a a_2$, centered at $a$ and belonging to a  $\Lambda $-factor 
of $A$.  
Let $e_i = aa_i$, and so $E(L) = \{e_1,e_2\}$.

Let $B$ be the graph-skeleton of the three-prism, say, 
$V(B) = \{1, \ldots , 6\}$ and $B$ is obtained from two disjoint triangles $123$ and $456$ by adding three 
new edges $14$, $25$, and $36$.

\begin{figure}
  \centering
   \includegraphics{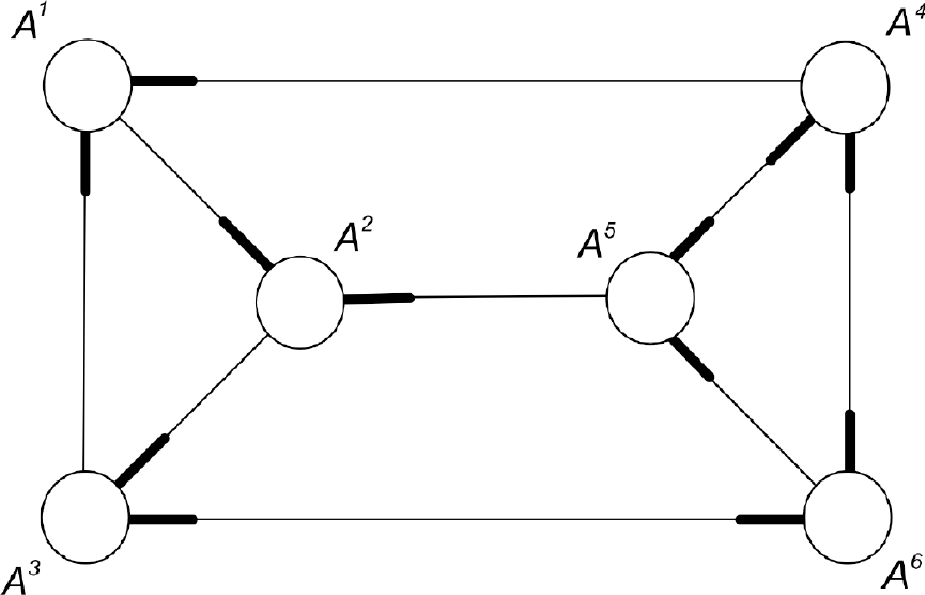}
  \caption{$(z2)$ $\Rightarrow$ $(z5)$}
  \label{fPrz2toz5}
\end{figure}

Let each $(A(v), a^v, a^v_1, a^v_1)$,  $v \in V(B)$ 
be a copy of $(A, a, a_1, a_2)$, 
and so edge $e^v_i = a^va^v_i$ in $A(v)$ is a copy 
of edge $e_i = aa_i$ in $A$, $i \in \{1,2\}$. 
We also assume that all $A(v)$'s are disjoint.
Let $G$ be a graph obtained from $B$ by replacing each 
$v \in V(B)$ by $A^v = A(v)- a^v)$
(see Fig. \ref{fPrz2toz5}).
Given $v \in V(B)$, let  $S(v) $ be the set of two edges $e'_i$ in $E'$ such that edge $e'_i$ is incident to vertex $a^v_i$ in $G$,  $\{i \in \{1,2\}$. 
We assume that each vertex $v$ in $B$ is replaced by 
$A^v$ (to obtain $G$) in such a way that
\\[0.5ex]
$S(x_1) = \{\alpha (13), \alpha (1\}$,
$S(x_2) = \{\alpha (21), \alpha (23)\}$,
$S(x_3) = \{\alpha (32), \alpha (36)\}$,
\\[0.5ex]
$S(y_1) = \{\alpha (45), \alpha (46)\}$,
$S(y_2) = \{\alpha (54), \alpha (56)\}$,
$S(y_3) = \{\alpha (63), \alpha (64)\}$.
\\[0.5ex]
In Figure \ref{fPrz2toz5} the edges in $S(v)$ are marked 
for every $v \in V(B)$.
\\[0.5ex]
\indent
By {\bf \ref{AasbB}}, $G$ is a cubic, 3--connected graph. 
Since $v(B) = 0\bmod 6$, clearly also $v(G) = 0\bmod 6$.
By $(z2)$, $G' = G - \alpha (36)$ has a $\Lambda $-factor, say $P$.

We know that 
$A$ has exactly one 3-vertex path 
$L = a_1 a a_2$ centered at $a$ and belonging to 
a  $\Lambda $-factor of $A$ and that each 
$(A^v, a^v, a^v_1, a^v_2)$ is a copy of $(A, a, a_1,a_2)$, 
and so $v(A^v - a^v) = 1 \bmod 6$.
Therefore by {\bf \ref{3cut}} $(a1)$ 
the  $\Lambda $-factor $P$ satisfies the following condition
for every $v \in V(B)$: 
\\[1ex]
${\bf c(v)}$ if $cmp (P^v) = 2$ then 
$v  \neg a$ and $v  \neg b$, where 
$\{\alpha (va), \alpha (vb)\} = S(v)$.
\\[1ex]
\indent
Obviously $|D^3| = |D^6| = 2$ in $G - \alpha (36)$.
Therefore $cmp (P^3) \le 2$ and $cmp (P^6) \le 2$.
Since $\alpha (36) \in S(3) \cap S(6)$, 
by conditions ${\bf c(3)}$ and ${\bf c(6)}$, 
$cmp (P^3) = cmp (P^6) = 1$.
Now by {\bf \ref{3cut}} $(a1)$, 
$x' \neg 3$ for some $x' \in \{1, 2\}$ and
$y' \neg 6$ for some $y' \in \{4, 5\}$.
\\[1ex]
${\bf (p1)}$ Suppose that $1\neg 3$.
Assume first that $\alpha (14) \not \in E(P)$.
Then $cmp(P^1) \le 2$. By {\bf \ref{3cut}} $(a1)$,
$cmp(P^1) = 2$. This contradicts ${\bf c(1)}$.
Thus we can assume that $\alpha (14) \in E(P)$.
\\[0.5ex]
${\bf (p1.1)}$ Suppose that $4 \neg 6$.

Suppose that $1\neg 4$. 
Then by {\bf \ref{3cut}} $(a1)$, $5 \neg 4$ and 
$5 \neg 2$. This contradicts ${\bf c(5)}$.

Now suppose that $4\neg 1$.
This contradicts ${\bf c(4)}$.
\\[0.5ex]
${\bf (p1.2)}$ Suppose that $5 \neg 6$.
Then $cmp (P^4) \le 2$.

Suppose that $1\neg 4$. 
By {\bf \ref{3cut}} $(a1.1)$, $cmp (P^4) = 1$.
Then $5 \neg 2$. This contradicts  ${\bf c(5)}$.

Now suppose that $4\neg 1$.
Then $cmp (P^4) = 2$. This contradicts ${\bf c(4)}$.
\\[1ex]
${\bf (p2)}$ Now suppose that $2\neg 3$.
Then $cmp(P^1) \le 2$.
By ${\bf c(2)}$ and {\bf \ref{3cut}} $(a1.3)$,
$x_1 \neg x_2$ (and $y_2 \neg x_2$).
Then by {\bf \ref{3cut}} $(a1.2)$, $cmp(P^1) = 2$.
This contradicts ${\bf c(1)}$.
\ep
\bs 
\label{z7toz5}
$(z7)$ $\Rightarrow$ $(z5)$.
\es

A proof of {\bf \ref{z7toz5}} can be obtained from 
the above {\bf Proof} of  {\bf \ref{z2toz5}} by using
$(z7)$ instead of $(z2)$ and by eliminating ${\bf (p.1.1)}$.
\bs 
\label{z1ottoz6}
$(z1)$ $\Leftrightarrow$ $(z6)$.
\es

\bp Obviously $(z1)$ $\Leftarrow$ $(z6)$ and 
$(z5) \Rightarrow$ $(z6)$. 
Now $(z1)$ $\Rightarrow$ $(z6)$ follows from 
{\bf \ref{z1ottoz2}} and {\bf \ref{z2toz5}}. 
\ep
\bs 
\label{z7toz8}
$(z7)$ $\Rightarrow$ $(z8)$.
\es

\noindent
{\bf Proof 1}.~
Suppose, on the contrary, that $(z7)$ is true but $(z8)$ 
is not true. Then there is a cubic 3-connected graph $A$
and a 3-path $L = a_1aa_2$ in $A$ such that 
$v(A) = 0\bmod 6$ and $A - L$ has no  
$\Lambda $-factor. Let $N(a, A) = \{a_1, a_2, a_3\}$.
Let $(A^i;  a^i, a^i_1, a_2, a^i_3)$, $i \in \{1,2\})$, be two copies of $(A; a, a_1, a_2, a_3)$ and $A^1$, $A^2$ 
be disjoint graphs, and so $L^i = a^i_1 a^i a^i_2$ in 
$A^i$ is a copy of $L = a_1 a a_2$ in $A$.
Let $G = A^1a^1\sigma a^2A^2$, where 
$\sigma : N(a^1, A^1) \to N(a^2, A^2)$ is a bijection such that $\sigma (a^1_i) = a^2_i$ for $i \in \{1, 2, 3\}$.
Let $H$ be the graph obtained from $G$ by subdividing
edge $a^1_j a^2_i$ by a new vertex $v_j$ for $j \in \{1,2\}$ 
and by adding a new edge $v_1v_2$ 
Obviously $G$ is a cubic 3-connected graph and
$v(G) = 0 \bmod 6$. By $(z7)$, $G - \{v_1v_2, a^1_3 a^2_3\}$ has a $\Lambda $-factor, say  $P$. 
Since $v(A^i - a^i) = -1\bmod 6$, clearly
$a^1_1 v_1 a^2_1$ and $a^1_2 v_1a^2_2$ are components of $P$. Then $A - \{a_1, a_2\} = A - L$ has a $\Lambda $-factor,
a contradiction.
\ep
\\

\noindent
{\bf Proof 2} (uses {\bf \ref{3cut}} $(a1)$).~
Suppose, on the contrary, that $(z7)$ is true but $(z8)$ 
is not true. Then there is a cubic 3-connected graph $A$
and a 3-path $L = a_1aa_2$ in $A$ such that 
$v(A) = 0\bmod 6$ and $A - L$ has no  $\Lambda $-factor. 
Let $N(a,A) = \{a_1,a_2,a_3\}$.

\begin{figure}
  \centering
   \includegraphics{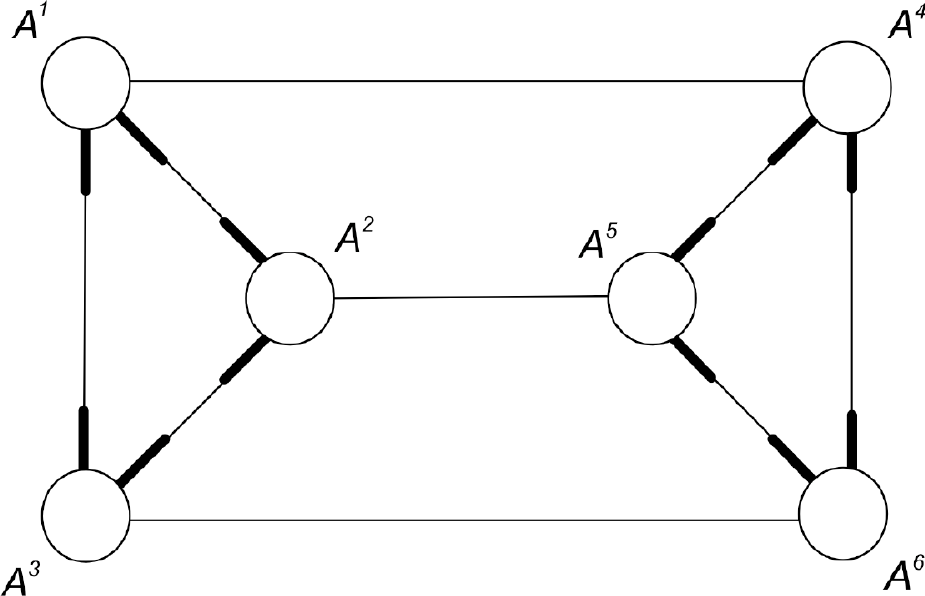}
  \caption{$(z7)$ $\Rightarrow$ $(z8)$}
  \label{fPrz7toz8}
\end{figure}

Let $B$,  $\{(A(v), a^v, a^v_1, a^v_2, a^v_3): v \in V(B)\}$, 
and $G$ be as in {\bf \ref{z2toz5}} (see Fig. \ref{fPrz7toz8}). 
 
Given $v \in V(B)$, let  $S(v) $ be the set of two edges $e'_i$ in $E'$ such that edge $e'_i$ is incident to vertex $a^v_i$ in $G$,  $\{1 \in \{1,2\}$. 
We assume that each vertex $v$ in $B$ is replaced by 
$A^v$ (to obtain $G$) in such a way that
\\[0.5ex]
$S(x_1) = \{\alpha (12), \alpha (13)\}$,
$S(x_2) = \{\alpha (21), \alpha (23)\}$,
$S(x_3) = \{\alpha (32), \alpha (31)\}$,
\\[0.5ex]
$S(y_1) = \{\alpha (45), \alpha (46)\}$,
$S(y_2) = \{\alpha (54), \alpha (56)\}$,
$S(y_3) = \{\alpha (64), \alpha (65)\}$.
\\[0.5ex]
In Figure \ref{fPrz7toz8}
the edges in $S(v)$ are marked 
for every $v \in V(B)$. 
\\[0.5ex]
\indent
Since $v(G) = 0\bmod 6$ and $G$ is cubic and 3-connected, by $(z7)$, $G - \{\alpha (14), \alpha (3,6)\}$ has a $\Lambda $-factor, say  $P$. 
By {\bf \ref{3cut}} $(a1)$, for every $v \in V(B)$, 
the  $\Lambda $-factor $P$ satisfies the following condition: 
\\[1ex]
${\bf c(v)}$ if $cmp (P^v) = 2$ then 
$v  \neg a$ and $v  \neg b$ where 
$\{\alpha (va), \alpha (vb)\} \ne S(v)$.
\\[1ex]
\indent
Obviously $|D^i| =  2$ in $G - \{\alpha (14), \alpha (3,6)\}$, 
and so  $cmp (P^i) \le 2$ for  $i \in \{1,3\}$.
Since $S(1) = \{12,13\}$ and $S(3) = \{31,32\}$, 
by conditions ${\bf c(1)}$ and ${\bf c(4)}$ we have: 
$cmp (P^1 = cmp (P^3) = 1$.
Now by {\bf \ref{3cut}} $(a1.1)$, $2 \neg 1$ and $2 \neg 3$.
This contradicts ${\bf c(1)}$.
\ep
\bs 
\label{z8ottoz7}
$(z8)$ $\Rightarrow$ $(z7)$.
\es
\bp
Obviously $(z8) \Rightarrow (z4)$.
By {\bf \ref{z4ottot2}}, $(z4)\Rightarrow (t2)$ and
by {\bf \ref{t2ottoz7}}, $(t2) \Rightarrow (z7)$.
\ep
\bs 
\label{t2ottoz7}
$(t2)$ $\Leftrightarrow$ $(z7)$.
\es

\bp (uses {\bf \ref{z1ottoz4}} and {\bf \ref{z4ottot2}}).
We first prove $(t2) \Rightarrow (z7)$.
Let $G$ be a cubic, 3-connected graph
with $v(G) = 0 \bmod 6$ and $a = a_1a_2$, $b = b_1b_2$ 
be two distinct edges of $G$.
Let $G'$ be the graph obtained from $G$ as follows:
subdivide edge $a_1a_2$ by a new vertex $a'$ and
edge $b_1b_2$ by a new vertex $b'$ and add a new edge 
$e = a'b'$. Then $G'$ is a cubic and 3-connected graph,
$v(G') = 2 \bmod 6$, and $G - \{a,b\} = G' - \{a',b'\}$.
By $(t2)$, $G' - \{a',b'\}$ has a $\Lambda $-factor.

Now we prove $(t2)\Leftarrow(z7)$.
Obviously $(z7) \Rightarrow (z1)$.
By {\bf \ref{z1ottoz4}},
$(z1) \Rightarrow (z4)$ and by {\bf \ref{z4ottot2}},
$(z4)\Rightarrow (t2)$.
Implication $(t2)\Leftarrow(z7)$ also follows from obvious
$(z8)\Rightarrow (z4)$, from $(z7)\Rightarrow (z8)$, 
(by {\bf \ref{z7toz8}}), and from $(z4)\Rightarrow (t2)$
(by {\bf \ref{z4ottot2}}).
\ep
\\[1.5ex]
\indent
Here is a direct proof of $(z7)\Rightarrow (t2)$.
\bs 
\label{z7tot2}
$(z7)\Rightarrow (t2)$.
\es

\bp Let $G$ be a cubic, 3-connected graph, 
$v(G) = 2\bmod 6$, $xy \in E(G)$, 
$N(x,G) = \{x_1,x_2, y\}$, and $N(y,G) = \{y_1,y_2, x\}$.
Let $G_1 = G - \{x, y\} \cup E_1$,
$G_2 = G - \{x, y\} \cup$, and
$G_3 = G - \{x, y\} \cup E_3$,
where $E_1 = \{x_1y_1,x_2y_2\}$,  
$E_2 = \{x_1y_2,x_2y_1\}$, and 
$E_3 = \{x_1x_2,y_1y_2\}$.
Obviously each $G_i$ is a cubic graph.
It is easy to see that since $G$ is 3-connected, 
there is $s \in \{1,2,3\}$ such that $G_s$ 
is 3--connected. Clearly $G_s - E_s= G - \{x, y\}$.
By $(z7)$, $G_s - E_s$ has a $\Lambda $-factor.
\ep
\bs 
\label{z1ottoz8}
$(z1)$ $\Leftrightarrow$ $(z8)$.
\es

\bp  Obviously $(z8)$ $\Rightarrow$ $(z1)$.
By {\bf \ref{z1ottoz4}},  $(z1)$ $\Rightarrow$ $(z4)$.
By {\bf \ref{z4ottot2}}, $(z4)$ $\Rightarrow$ $(t2)$.
By {\bf \ref{t2ottoz7}}, $(t2)$ $\Rightarrow$ $(z7)$.
By  {\bf \ref{z7toz8}}, $(z7)$ $\Rightarrow$ $(z8)$.
Therefore $(z1)$ $\Rightarrow$ $(z8)$.
\ep
\bs 
\label{z8tof1}
$(z8)$ $\Rightarrow$ $(f1)$.
\es

\bp 
Let $G$ be a cubic 3-connected graph, 
$v(G) = 4 \bmod 6$, $x \in V(G)$, and 
$N(x,G) = \{x_1,x_2, x_3\}$.
Let $G'$ be the graph obtained from $G$ by replacing $x$ 
by a triangle $T$ with $V(T) =  \{x'_1, x'_2, x'_3\}$ so that $x_ix'_i \in E(G')$, $i \in \{1,2,3\}$.
Since $v(G) = 4 \bmod 6$, clearly $v(G') = 0 \bmod 6$.
Consider the 3-vertex path  $L' = x'_1x'_2x'_3$ in $G'$. 
By $(z8)$, $G' - L'$ has a $\Lambda $-factor, say  $P'$. Obviously $P' - L'$ is 
a $\Lambda $-factor of $G' - L'$ and $G - x = G' - L'$.
\ep

\bs 
\label{z8tof2}
$(z8)$ $\Rightarrow$ $(f2)$.
\es

\begin{figure}
  \centering
   \includegraphics{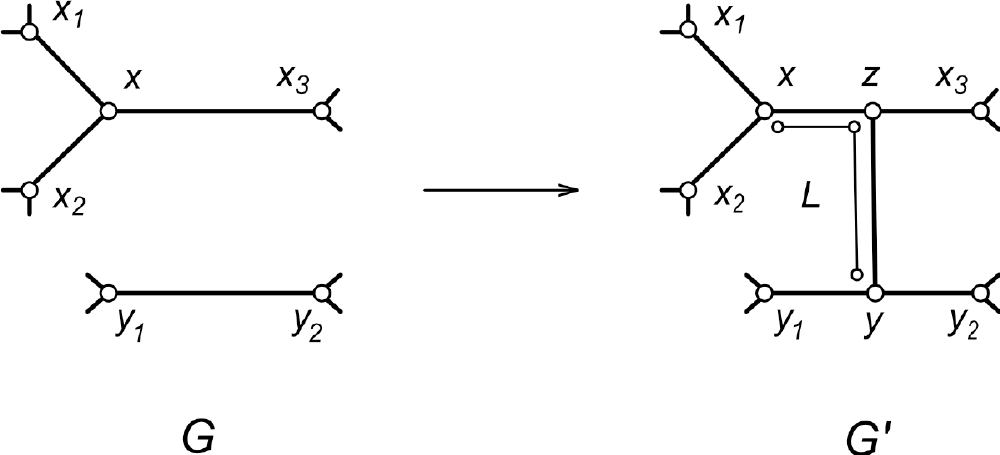}
  \caption{$(z8)$ $\Rightarrow$ $(f2)$}
  \label{fPrz8tof2}
\end{figure}

\bp (uses {\bf \ref{z8tof1}}).
Let $G$ be a cubic 3-connected graph, $v(G) = 4 \bmod 6$,
$x \in V(G)$, and $e = y_1y_2 \in E(G)$. We want to prove that if $(z8)$ is true, then $G - \{x,e\}$ has a 
$\Lambda $-factor. 
If $x  \in \{y_1,y_2\}$, then $G - \{x,e\} = G - x$, and therefore by {\bf \ref{z8tof1}}, our claim is true. 
So we assume that $x  \in \{y_1,y_2\}$.
Let $N(x,G) = \{x_1,x_2,x_3\}$.
Let $G'$ be the graph obtained from $G$ by subdividing
edge $y_1y_2$ by a vertex $y$ and edge $xx_3$ by a vertex
$z$ and by adding a new edge $yz$
(see Fig. \ref{fPrz8tof2}).
Since  $v(G) = 4 \bmod 6$, clearly $v(G') = 0 \bmod 6$.
Since $x \not \in \{y_1,y_2\}$ and $G$ is cubic and 3-connected, $G'$ is also  cubic and 3-connected. Obviously
$L = xzy$ is a 3-vertex path in $G'$ and $G - \{x,e\} =
G - \{x, y_1y_2\} = G' - L$. 
By $(z8)$, $G' - L$ has a $\Lambda $-factor.
\ep
\bs 
\label{f2tot4}
$(f2)$ $\Rightarrow$ $(t4)$.
\es

\bp Let $G$ be a cubic, 3-connected graph, 
$v(G) = 2\bmod 6$, $x \in V(G)$, and 
$N(x,G) = \{x_1,x_2, x_3\}$.
Let $G'$ be the graph obtained from $G$ by replacing $x$ 
by a triangle $T$ with $V(T ) =  \{x'_1, x'_2, x'_3\}$ so that $x_ix'_i \in E(G')$, $i \in \{1,2,3\}$.
Since $v(G) = 2 \bmod 6$, clearly $v(G') = 4 \bmod 6$.
By $(f2)$, $G' - \{x'_i,x'_jx'_k\}$ has a $\Lambda $-factor, 
say  $P_i$ where $\{i, j, k\} = \{1, 2, 3\}$. 
Since $x_jx'_j$ and $x_kx'_k$ are dangling edges in 
$G' - \{x'_i,x'_jx'_k\}$, clearly $x_jx'_j, x_kx'_k \in E(P_i$ and 
$d(x'_j, P_i) = d(x'_k, P_i) =1$. Let $L_j$ and $L_k$ be the components of $P_i$ containing $x_jx'_j$ and $x_kx'_k$, respectively. Then $E(L_j) \cup E(L_k)$ induces in $G$ 
a 5-vertex path $W_i$ such that  $x$ is the center vertex of $W_i$ and $x_ix'_i \not \in E(W_i)$.
\ep
\bs 
\label{t3toz1}
$(t3)$ $\Rightarrow$ $(z1)$.
\es

\bp ({\bf \ref{3cut}} $(a1.3)$).
 Let $G = Y(A^1,a^1; A^2,a_2; A^3,a^3)$
(see Fig. \ref{fY}).
Suppose, on the contrary, that $(t3)$ is true but $(z1)$ 
is not true, i.e. there is a cubic 3-connected graph $A$  
such that $v(A) = 0 \bmod 6$ and $A$ has no 
$\Lambda $-factor. Let $a \in V(A)$.
Let $(A^i,a^i)$ above be a copy of $(A, a)$ for $i \in \{1,2\}$
and let $(A^3, a^3)$ be such that $v(A_3) = 2 \bmod 6$.
Obviously $v(G) = 2 \bmod 6$.
By $(t3)$, $G$ has a 5-vertex path $W$ such that $z_3$ is the center vertex of $W$ and $G - W$ has a 
$\Lambda $-factor, say $P$.
Obviously $(A^i - a^i) \cap W = \emptyset $ for some 
$i \in \{1,2,3\}$.

Suppose that $(A^3 - a^3) \cap W = \emptyset $.
Then $W$ has an end-edge in $A^1 - a^1$ and 
in $A^2 - a^2$.
Since $A$ has no $\Lambda $-factor, 
by {\bf \ref{3cut}} $(a1.3)$, $D^i - e(W) \subseteq E(P)$ for 
$i \in \{1,2\}$. Then $P$ is not a $\Lambda $-factor of 
$G - W$, a contradiction.

Now suppose that $(A^3 - a^3) \cap W \ne \emptyset $.
By symmetry, we can assume 
that $(A^2 - a^2) \cap W = \emptyset $.
Then $W$ has an end--edge in $A^1 - a^1$ and 
in $A^3 - a^3$.
Then by {\bf \ref{3cut}} $(a1.3)$,  
$Cmp(P^1) = \{L_1,L_2\}$, where $L_1$  has an end-vertex 
in  $A^1 - a^1$ and $L_2$ has an end-edge in $A^1 - a^1$.
By symmetry, we can assume that $a^1_iz_i \in E(L_i)$ 
for $i \in \{1,2\}$. 
Then $L_1 = a^1_1z_1y$, where $y \in \{a^2_1,a^3_1\}$ and
$z_2$ is of degree one in $P$. 
Then $P$ is not a $\Lambda $-factor of 
$G - W$, a contradiction.
\ep
\bs 
\label{z8tof4}
$(z8)$ $\Rightarrow$ $(f4)$.
\es

\bp Let $G$ be a cubic, 3-connected graph, 
$v(G) = 4 \bmod 6$, $x \in V(G)$, and 
$N(x,G) = \{x_1,x_2, x_3\}$.
Let $G'$ be the graph obtained from $G$ by replacing $x$ 
by a triangle $ \Delta $ with $V(\Delta ) =  \{x'_1, x'_2, x'_3\}$ so that $x_ix'_i \in E(G')$, $i \in \{1,2,3\}$.
Since $v(G) = 4 \bmod 6$, clearly $v(G') = 0 \bmod 6$.
Consider the 3-vertex path  $L'_i = x_jx'_jx'_k$ in $G'$, 
where $\{i, j, k\} = \{1, 2, 3\}$. 
By $(z8)$, $G' - L'_i$ has a $\Lambda $-factor, say  $P_i$. Since $x_ix'_i$ is a dangling edge in 
$G' - L'_i$, clearly $x_ix'_i \in E(P_i$ and $d(x'_i, P_i) = 1$. 
Let $L_i$ be the components of $P_i$ containing $x_ix'_i$. Then $E(L_i) \cup x_jx'_j$ induces in $G$ a 4-vertex path $Z_k$ such that $x_i$ is an inner vertex of $Z_k$ and $x_kx'_k \not \in E(Z_k)$.
\ep
\\[1.5ex]
\indent
Let $H'$ be a tree such that 
$V(H') = \{x,y\} \cup (b^j: j \in \{1,2,3,4\}$ and 
$E(G) = \{xy, b^1x, b^2x, b^3y, b^4y\}$.
Let $H_i$, $i \in \{1,2,3\}$ be three disjoint copies of $H'$ 
with $V(H_i) = \{x_i,y_i\} \cup (b^j_i: j \in \{1,2,3,4\}$.
Let $H$ be obtained from these three copies by identifying 
for every $j$ three vertices $b^j_1$, $b^j_2$, $b^j_2$ with a new vertex $z^j$.
Let $A^i$, $i \in \{1,2,3,4\}$, be a cubic 
graph, $a^i \in V(A^i)$ and let
$G = H(A^1,a^1; A^2,a_2; A^3,a^3; A^4,a^4)$ be the graph obtained from $H$ by replacing each $z^j$ by $A^j - a^j$ assuming that all $A^i$'s are disjoint, 
\bs 
\label{f3toz1}
$(f3)$ $\Rightarrow$ $(z1)$. 
\es

\bp (uses {\bf \ref{A(B)}}).
Let $G = H(A^1,a^1; A^2,a_2; A^3,a^3; A^4,a^4)$,
where each $A^i$ is a cubic 3-connected graph.
Since $H$ is cubic and 3-connected, by {\bf \ref{A(B)}},
$G$ is also cubic and 3-connected.

Suppose, on the contrary, that $(f3)$ is true but $(z1)$ 
is not true, i.e. there is a cubic 3-connected graph $A$
such that $v(A) = 0 \bmod 6$ and $A$ has no 
$\Lambda $-factor. Let $a \in V(A)$.
Let $(A^i,a^i)$ be a copy of $(A,a)$ for $i \in \{1,2,3\}$ and 
let  $v(A^4) = 2 \bmod 6$.
Obviously $v(G) = 4 \bmod 6$.
It is easy to see that $G - Z$ has no $\Lambda $-factor
for every 4-vertex path $Z$ in $G$ such that
$y_1$ is an inner vertex of $Z$. 
This contradicts $(f3)$.
\ep
\\[1.5ex]
\indent
Obviously $(f4)$ $\Rightarrow$ $(f3)$. Therefore from 
{\bf \ref{f3toz1}} we have: $(f4)$ $\Rightarrow$ $(z1)$.
Below we give a direct proof of this implication.
\bs 
\label{f4toz1}
$(f4)$ $\Rightarrow$ $(z1)$. 
\es

\bp 
Let $G = Y(A^1,a^1; A^2,a_2; A^3,a^3)$
(see Fig. \ref{fY}).
Suppose, on the contrary, that $(f4)$ is true but $(z1)$ 
is not true, i.e. there is a cubic 3-connected graph $A$
such that $v(A) = 0 \bmod 6$ and $A$ has no 
$\Lambda $-factor. Let $a \in V(A)$.
Let $(A^3,a^3)$ above be a copy of $(A, a)$
and let $(A^i, a^i)$ for $i \in \{1,2\}$ be copies of 
$(B,b)$ where $B$ is a  cubic 3-connected graph,
$v(B) = 2 \bmod 6$, and $b \in V(B)$.
Obviously $v(G) = 4 \bmod 6$.
By $(f4)$, $G$ has a 4-vertex path $Z$ such that $z_3$ 
is an inner  vertex of $Z$,  
$(A^3 - a^3) \cap Z = \emptyset $, and $G - Z$ has 
a $\Lambda $-factor, say $P$.
Since $A^3$ has no $\Lambda $-factor,
clearly $P$ is not a $\Lambda $-factor of $G - Z$, 
a contradiction.
\ep
\\[1.5ex]
\indent
Implication $(z8) \Rightarrow (f1)$ follows from obvious
$(z8) \Rightarrow$ $(z1)$ and from $(z1) \Rightarrow (f1)$,
by {\bf \ref{z1ottof1}}. 
It also follows from obvious
$(f2) \Rightarrow$ $(f1)$ and from $(z8) \Rightarrow (f2)$,
by {\bf \ref{z8tof2}}. Below we give a direct proof of this implication.
\bs 
\label{z1ottof6}
$(f6)$ $\Rightarrow$ $(f5)$ $\Rightarrow$ $(f4)$ $\Rightarrow$ $(z1)$ $\Rightarrow$ $(f6)$. 
\es

\bp (uses {\bf \ref{f4toz1}} and {\bf \ref{z1ottoz8}}).
Obviously $(f6)$ $\Rightarrow$ $(f5)$ $\Rightarrow$ $(f4)$. By {\bf \ref{f4toz1}}, $(f4)$ $\Rightarrow$ $(z1)$. 
Therefore $(f6)$ $\Rightarrow$ $(z1)$. It remains to prove
$(z1)$ $\Rightarrow$ $(f6)$. 
By {\bf \ref{z1ottoz8}}, $(z1)$ $\Rightarrow$ $(z8)$.
Thus it is sufficient to show that $(z8)$ $\Rightarrow$ $(f6)$.
Let $G$ be a cubic 3-connected graph, 
$v(G) = 4 \bmod 6$, and $xyz$ is a 3-vertex path in $G$.
Let $N(y,G) = \{x,z,s\}$ and $G'$ be obtained from $G$ by
subdividing edges $yz$ and $ys$ by new vertices $z'$ and $s'$, respectively, and by adding a new edge $s'z'$.
Then $G'$ is a cubic 3-connected graph and $v(G') = 0 \bmod 6$. Consider the 3-vertex path  $L' = zz's'$ in $G'$. 
By  $(z8)$, $G'$ has a $\Lambda $-factor $P'$ containing $L'$. Since vertex $y$ has degree one in $G' - L'$, clearly
$P'$ has a 3-vertex path $Q' = yxq$. 
Let $Z$ be the 4-vertex path $qxyz$ in $G$.
Then $P' - (P' \cup Q')$ is a $\Lambda $-factor in $G - Z$.
\ep
\bs 
\label{z9ottoz1}
$(z9)$ $\Leftrightarrow$ $(z1)$.
\es

\bp 
Obviously $(z9)$ $\Rightarrow$ $(z8)$.
By the above claims,  $(z1)$, $(z8)$, $(t4)$, and $(f6)$ are equivalent. So we can use these claims to prove $(z9)$.
Let $G$ be a cubic 3-connected graph, $K$ a 3-edge cut of $G$, $S \subset K$ and $|S| = 2$, and $v(G) = 0 \bmod 6$. 
If the edges of $K$ are incident to the same vertex $x$ in $G$ (i.e. $K = D(x,G)$), then by $(z8)$, $G$ has  a $\Lambda $-factor $P$ of $G$ such that $E(P) \cap K = S$, and so our claim is true.
So we assume that the edges in $K$ are not incident to the same vertex in $G$. 
Then since $G$ is cubic and 3-connected, clearly 3-edge cut  $K$ is matching. 
Let $A$ and $B$ be the two component of $G - K$. By the above arguments, we assume that 
$v(A)\ne 1$ and $v(B) \ne 1$.
Let $A^b$ be the graph obtained from $G$ by identifying the vertices of $B$ with a new vertex $b$ and similarly,
$B^a$ be the graph obtained from  $G$ by identifying the vertices of $A$ with a new vertex $a$, and so $A = A^b - b$ and $B = B_a - a$. Then $D(b, A_b) = D(a,B_a) = D$, and
$S \subset D$. 
Let $S = \{a_1b_1, a_2b_2\}$, where 
$\{a_1,a_2\} \subset V(A)$ and $\{b_1,b_2\} \subset V(B)$.
Obviously  $S$ forms a 3-vertex path $S_A = a_1ba_2$  in $A^b$ and a 3-vertex path $S_B = b_1ab_2$ in $B^a$. 
Since $G$ is cubic and 3-connected, both $A^b$ and $B^a$ are also cubic and 3-connected.
Since $v(G) = 0 \bmod 6$, there are two possibilities
(up to symmetry):
\\[0.5ex]
$(c1)$ $v(A^b) = 0 \bmod 6$ and $v(B^a) = 2 \bmod 6$ and 
\\[0.5ex]
$(c2)$ $v(A^b) = v(B^a) = 4 \bmod 6$.
\\[1ex]
\indent
Consider case $(c1)$. By $(z8)$, $A^b - S_A$ has a $\Lambda $-factor $P_A$. By $(t4)$, 
$B^a$ has a 5-vertex path $W$ such $a$ is the center vertex of $W$, $S_B \subset W$, and 
$B^a - W$ has a $\Lambda $-factor $P_B$.
Then $E(P_A) \cup E(P_B) \cup S$ induces 
a $\Lambda $-factor $P$ in $G$ such that 
$E(P) \cap K = S$.
\\[1ex]
\indent
Now consider case $(c2)$. 
By $(f6)$, we have:
\\[0.5ex]
$(a)$ 
$A^b$ has a 4-vertex path $Z_A$ such that $a_1$ is an end-vertex of $Z_A$, $S_A \subset Z_A$, and $A^b - Z_A $
has a $\Lambda $-factor, say $P_A$, and similarly,
\\[0.5ex]
$(b)$
$B^a$ has a 4-vertex path $Z_B$ such that $b_2$ is an end-vertex of $Z_B$, $S_B \subset Z_B$, and $B^a - Z_B$ has a $\Lambda $-factor, say $P_B$.
\\[0.7ex]
\indent
Then $E(P_A) \cup E(P_B) \cup S$ induces 
a $\Lambda $-factor $P$ in $G$ such that 
$E(P) \cap K = S$.
\ep

\section{On almost cubic graphs with no $\Lambda $-factors}
\label{AlmostCubic}
\indent

In the previous section we indicated that some claims in 
{\bf \ref{3-con}} (equivalent to $(z1)$) are best possible in some sense. In this section we describe constructions that provide some additional  facts of this nature.

Let $G = Y(A^1,a^1; A^2,a^2; A^3,a^3)$
(see Fig. \ref{fY}), where $A^3$ is the graph having two vertices 
$x$,  $a^3$ and three parallel edges with the end-vertices
$x$,  $a^3$, and so $A^3 - a^3 = x$.
\\[0.5ex]
\indent
If $v(A^1) = v(A^2) = 0 \bmod 6$, then
\\[0.5ex]
$(a1)$
$v(G) = 2\bmod 6$ and $G - (N(x,G)\cup x \cup y)$
has no $\Lambda $-factor for every vertex $y$ in 
$G - (N(x,G)\cup x)$ adjacent to a vertex in $N(x,G)$
$($see also {\bf \ref{3cut}} $(a1.3)$ and Fig \ref{A1}$)$.
\\[0.5ex]
\indent
If $v(A^1) = 2$ and $v(A^2) = 4 \bmod 6$, then 
\\[0.5ex]
$(a2)$
$v(G) = 4\bmod 6$ and $G - (x \cup N(x))$ has no 
$\Lambda $-factor
$($see also {\bf \ref{3cut}} $(a2.3)$ and Fig \ref{A2}$)$.
\\[1ex]
\indent 
Thus from the above construction we have:
\bs 
\label{G,x,W}
There are infinitely many pairs $(G,x)$ such that
$G$ is a cubic 3-connected graph, $x \in V(G)$, and
 $(G,x)$ satisfies $(aj)$ above, $j \in \{1,2\}$.
 \es 
 
Using {\bf \ref{G,x,W}}, one can also prove  the following.
\bs 
\label{G,P,K}
There are infinitely many cubic 3-connected graphs $G$ 
such that $v(G) = 0 \bmod 6$ and
$|E(P) \cap K| \in \{1,2\}$
for every $\Lambda $-factor $P$ of $G$ and every  
3-edge cut $K$ of $G$.  
\es

Now we want to define the class ${\cal F}$ of graphs $G$ that have some special $\Lambda $-packing properties and that
are `almost' cubic.
 Using these graphs we will construct cubic 3-connected graphs mentioned in our above remark $(r1)$
concerning the result in {\bf \ref{3-con}}.

Let $L(G)$ denote the set of leaves (i.e. of vertices having degree one) of  a graph $G$. If $T$ is a subgraph  of $G$, then let $N(T,G)$ be the set of vertices in $G$ adjacent to some vertices in $T$ and, as above,  $D(T,G)$ the set of edges in $G$ having exactly one end in $T$.

\begin{figure}
  \centering
   \includegraphics{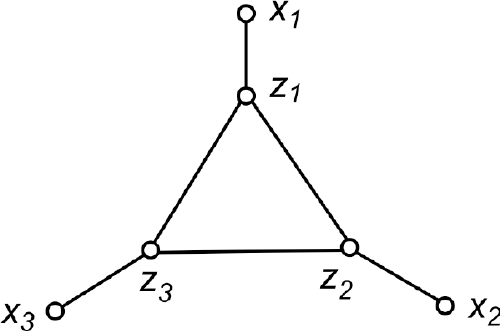}
  \caption{Graph $Y$}
  \label{fY0}
\end{figure}

\begin{figure}
  \centering
  \includegraphics{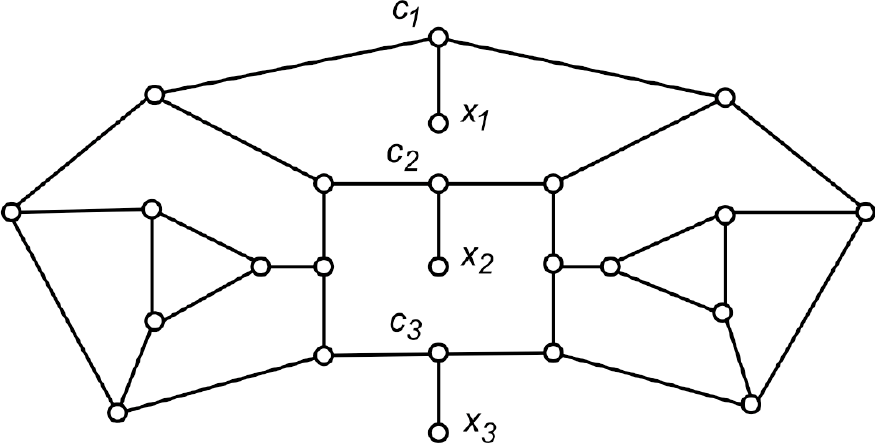}
  \caption{Graph $Z$}
  \label{fY1}
\end{figure}

First we define two special graphs $Y$ and $Z$.
Let $Y$
be the graph obtained from a triangle $T$ with $V(T) = \{z_1, z_2, z_3\}$ by adding three new vertices 
$x_1$, $x_2$, $x_3$ and three new edges 
$x_1y_1$, $x_2y_2$, $x_3y_3$), and so 
$x_1$, $x_2$, $x_3$ are the leaves of $Y_0$ 
(see Fig \ref{fY0}.
Let $Y'$ be a copy of $Y$
with the leaves $y_1$, $y_2$, $y_3$. Let $S$ be the graph obtained from $Y'$ by adding  six new vertices $\{s_i, s'_i: i \in \{1,2,3\}\}$ and nine new edges
$\{s'_iy_j, s'_iy_k: \{i,j,k\} =  \{1,2,3\}\} \cup \{s_is'_i:  i \in \{1,2,3\}\}$, and so $s_1$, $s_2$, $s_3$ are the leaves of $S$.
Let $A$ and $B$ be two disjoint copies of $S$ with the leaves
$a_1$, $a_2$, $a_3$ and $b_1$, $b_2$, $b_3$, respectively.
Let $Z$
be the graph obtained from $A$ and $B$ by identifying $a_i$ and $b_i$ with a new vertex $c_i$, $i \in \{1,2,3\}$,  and by adding three new vertices $x_1$, $x_2$, $x_3$ and three new edges $x_1c_1$, $x_2c_2$, $x_3c_3$, and so $x_1$, $x_2$, $x_3$ are the leaves of $Z$
(see Fig. \ref{fY1}).

Now we are ready to define the class of graphs ${\cal F}$ recursively.
First we assume that $Y$
and $Z$
are in ${\cal F}$.
Suppose that $A$ and $B$ are disjoint graphs such that
$A$ has a triangle $T$, $N(T,A) = \{a_1,a_2,a_3\}$, and
$L(B) = \{b_1,b_2,b_3\}$. 
Let $A(T,B)$ be a graph obtained from $A$ by {\em replacing} its triangle $T$ by $B - L(B)$, i.e. 
$A(T,B)$ is obtained from $A - T$ and $B$  by identifying 
each $a_i$ with $b_i$, $i \in \{1,2,3\}$.
Now if $A, B \in {\cal F}$, then we assume that also
$A(T,B) \in {\cal F}$.
\\[1.5ex]
\indent
It is easy to prove 
that the graphs in ${\cal F}$ have the following simple properties.
\bs
\label{CalYpropty}
Let $G \in {\cal F}$ and $G \ne Y$. Then
\\[0.5ex]
$(a1)$ $|L(G)| = 3$ and if $x \in V(G - L(G))$, then 
$d(x,G) = 3$,
\\[0.5ex]
$(a2)$ $G$ has triangles and if $T$ is a triangle of $G$, then
\\[0.5ex]
\indent
$(a2.1)$
$|N(T,G)| = |D(T,G) = 3$, 
\\[0.5ex]
\indent
$(a2.2)$
$N(T,G)$ induces in $G$ the subgraph with no edges and
$D(T,G)$ is a 3-edge cut-matching  in $G$,
\\[0.5ex]
\indent
$(a2.3)$ there is a unique 6-cycle $C$ in $G$ such that $N(T,G) \subset V(C)$ and $D(T\cup C,G)$ is a 
3-edge cut-matching in $G$ {\em (we put  $C = C(T,G)$ and
$D(T\cup C,G) = M(T,G)$)}. 
\es

If $F \in {\cal F}$, then let 
$\dot{F}$ denote the graph obtained from $F$ by identifying the three leaves with a new vertex $x$,   
$\bar{F}$
the graph obtained from $F$ by adding 
the triangle $T$ with the vertex set $L(F)$
 and
$\ddot{F}$ the graph obtained from $\bar{F}$ by adding 
a new vertex $z$, by subdividing every edge $e$ in $T$ 
with a new vertex $v_e$, and by adding three new edges 
$zv_e$, $e \in E(T)$.
\\[1.5ex]
\indent
It is easy to see the following.
\bs 
\label{dotYcubic3con} 
Let $F \in {\cal F}$. Then
$\dot{F}$, $\bar{F}$, and $\ddot{F}$ are cubic 3-connected graphs.
\es 

Now we can describe some $\Lambda $-packing properties of $F$, $\dot{F}$, $\bar{F}$, and  $\ddot{F}$ for $F \in {\cal F}$.
\bs 
\label{CalY} 
Let $F \in {\cal F}$. Then
\\[1ex]
$(a1)$
$v(F) = 0\bmod 6$ and
$F$ has no $\Lambda $-factor, 
\\[1ex]
$(a2)$ 
$v(\dot{F}) = 4 \bmod 6$ and
$\dot{F} - (N(x,\dot{F}) \cup x \cup X)$  has no 
$\Lambda $-factor for every $X \subset V(G)$ such that 
$|X| = 3$ and $X$ is matched with $N(x)$ in  $\dot{F}$
{\em (see also {\bf \ref{3cut}} $(a2)$ and Fig. \ref{A2})},
\\[1ex]
$(a3)$ 
$v(\bar{F}) = 0\bmod 6$ and
$\bar{F} - E(T)$ has no $\Lambda $-factor, where $T$ is the triangle in $\bar{F}$ with $V(T) = L(F)$, and
\\[1ex]
$(a4)$ 
$v(\ddot{F}) = 4 \bmod 6$ and
$\ddot{F} - (N(z,\ddot{F}) \cup z)$ has no $\Lambda $-factor.
\es

\bp (uses {\bf \ref{3cut}} $(a2)$ and {\bf \ref{CalYpropty}}).
Claims $(a2)$, $(a3)$, and $(a4)$ follow from $(a1)$.
We prove calim $(a1)$. Obviously, $v(F) = 0\bmod 6$ and our claim is obviously true for $Y$ and $Z$. Suppose, on the contrary, that $(a1)$ is not true.
Let $G$ be a vertex minimum counterexample, and so
$G \in {\cal F}$ and $G$ has a $\Lambda $-factor, say $P$.  
By definition of ${\cal F}$, we have: $G = A(T,B)$ for some 
$A, B \in {\cal F}$ and a triangle $T$ in $A$.
By  {\bf \ref{CalYpropty}}, there exist $M = M(T,A)$ and 
$C = C(T,A)$.
Obviously, $v(B \cup S)) = 0 \bmod 3$. 
Therefore 
$(P, M)$ satisfies one of the conditions  in {\bf \ref{3cut}} $(a2)$ (see Fig. \ref{A2}).
Let $Q = P \cup S$ and $B' = (B \cup C) - Q$. Then 
$P_1 = P \cap B'$  is a $\Lambda $-factor in $B'$ and
$P_2 = P - P_1$ is a  $\Lambda $-factor in $G - B'$.

Suppose that $(P, M)$ satisfies 
conditions 
$(a2.2)$ with  $E(S) \cap E(P) \ne \emptyset $
(and so $|E(C) \cap E(P)| \in \{1,3\}$) or $(a2.1)$ or $(a2.2)$.  
Then $T' = (T \cup C \cup D(T,A)) - Q$ has a $\Lambda $-factor $P'_1$. Therefore $P'_1 \cup P_2$ is a $\Lambda $-factor in $A$. However, $A \in {\cal F}$ and $v(A) < v(G)$.
Therefore the counterexample $G$ is not vertex minimum, a contradiction. 

\begin{figure}
  \centering
  \includegraphics{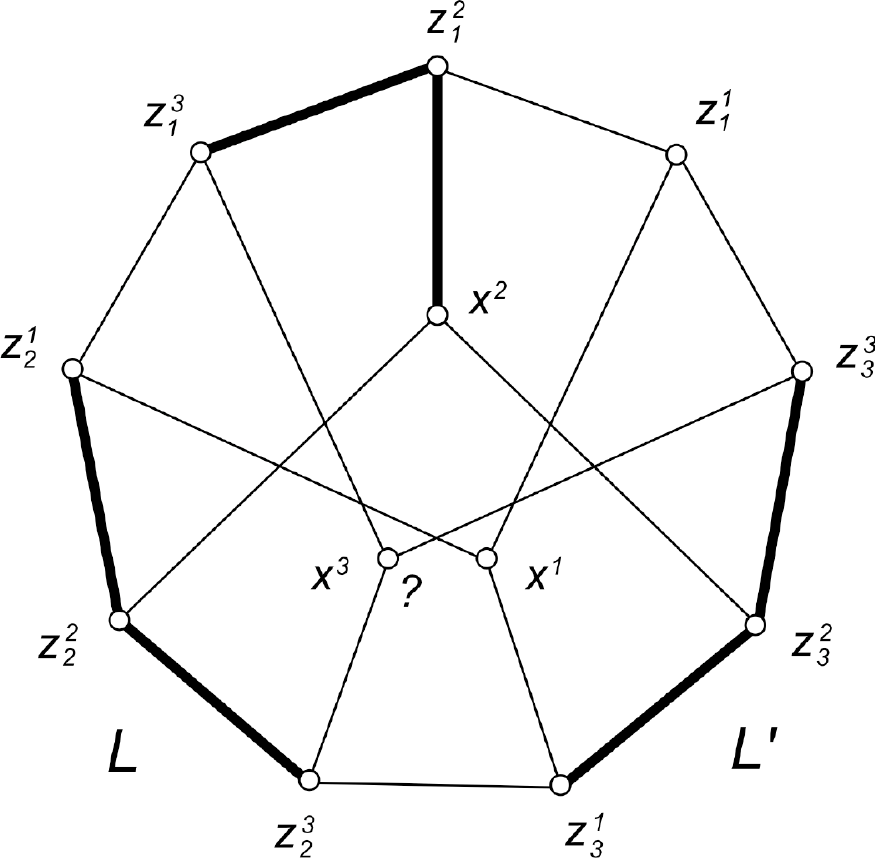}
    \caption{$(R_1, L, L')$, $v(R_1) = 12$}
  \label{fR1}
\end{figure}

Now suppose that $(P, M)$ satisfies 
condition $(a2.2)$ with  $E(C) \cap E(P) = \emptyset $.
Then $B' = B$, and so $P_1$ is a $\Lambda $-factor in $B$.
However, $B \in {\cal F}$ and $v(B) < v(G)$.
Therefore again the counterexample $G$ is not vertex minimum, a contradiction. 
\ep
\\[1.5ex]
\indent
Now we describe a sequence (mentioned in the above remark $(r2)$) of cyclically 6-connected graphs $G$ with two disjoint 3-vertex paths $L$, $L'$ such that
$v(G) = 0 \bmod 6$ and $G - (L\cup L')$ has no $\Lambda $-factor.
Let $C_s$ be a cycle with $9s$ vertices, $s \ge 1$ and let
$\{L_k: k \in \{1, \ldots , 3s\}$ 
be a $\Lambda $-factor of $C_s$, where $L_i = (z_i^1z_i^2 z_i^3)$.
Let $R_s$ be the graph obtained from $C_s$ by adding the set  
$\{x^j_i:  i \in \{1, \ldots , s\}, j \in \{1,2,3\}\}$ 
of $3s$ new vertices and  the set 
$\{x^j_iz^j_i, x^j_iz^j_{i+s}, x^j_iz^j_{i+2s}:  i \in \{1, \ldots , s\}, 
j \in \{1,2,3\}\}$ 
of $9s$ new edges 
(see, for example,  $(R_1, L, L')$ in Fig. \ref{fR1}). 
\\[1.5ex]
\indent
It is easy to prove the following
\bs
\label{G,L,L'} 
Let $R_s$ be the graph described above, $s \ge 1$, and $\{L,L'\} \subset \{L_i,L_{i+s}, L_{i+2s}\}$ for some 
$i \in \{1, \ldots , s\}$. Then
\\[0.5ex]
\indent
$(a1)$
$R_1$ is a cubic cyclically 5-connected graph, $R_s$ is a cubic cyclically 6-connected graph for $s \ge 2$,  
$v(R_s) = 12 s$, and 
\\[0.5ex]
\indent
$(a2)$
$R_s - (L \cup L')$ has no $\Lambda $-factor.
\es 

Using operation
$Aa\sigma bB$ 
(see Fig. \ref{fAasbB}), {\bf \ref{3cut}}, and  
{\bf \ref{CalY}}, it is easy to prove the following.
\bs 
\label{z9ottoz1}
There are infinitely many pairs $(G,K)$ such that $G$ is 
a cubic 3-connected graph,
$K$ is a 3-edge cut of $G$, $v(G) = 0 \bmod 6$, and 
$|E(P) \cap K| \not \in \{0,1\}$ for every 
$\Lambda $-factor $P$ of $G$.
\es

\section{On a  $\Lambda $-factor homomorphism in cubic graphs}
\label{homomorphism}

\indent

Let, as in Section \ref{constructions}, 
$G = B\{(A(u),a^u): u \in V(B)\}$,
$N(a^u,A(u)) = N^u = \{a^u_1,a^u_2,a^u_3\}$, 
$A^u = A(u) - a^u$,
and $E' = E(G) \setminus \cup \{E(A^v):  v \in V(B)\}$.
As we mentioned above, there is a unique bijection
$\alpha : E(B) \to E'$ such that if $uv \in E(B)$, then 
$\alpha (uv)$ is an edge in $G$ having one end-vertex in
$A^u$ and the other end-vertex in $A^v$.

Let $D(A^u,G) = D^u = \{a^u_ib^u_i: i \in \{1,2,3\}\}$,
and so $D^u$ forms a matching in $G$. 
Then $A^u \cup D^u$ is a subgraph of $G$.
Let $A'(u)$ be the graph obtained from $A^u \cup D^u$ 
by adding the triangle $T^u$ 
with $V(T^u = \{b^u_1, b^u_2, b^u_3\}$.
Suppose that $P$ is a $\Lambda $-factor of $B$ and 
$vuw$ is a 3-vertex path in $P$.
We need the following additional notation:
\\[0.7ex]
$V_s(P)$ is the set of vertices of degree $s$ in $P$ (and 
so $s \in \{1,2\}$),
\\[0.7ex]
$A_1(v,P) = A^v - End(\alpha (uv))$ and  
$A_2(u,P) = A'(u) - b$, where $b$ is the vertex in $V(T^u)$ 
that is incident to no edge in $\{\alpha (us), \alpha (uv)\}$,
\\[0.7ex]
$\Gamma (G, P)$ is the set of  
$\Lambda $-factors $Q$ of $G$ such that 
$E(P) = \{\alpha ^{-1}(e): e \in E'(G) \cap E(Q)\}$, 
\\[0.7ex]
$\Gamma (H)$ is the set of $\Lambda $-factors of a graph $H$, and 
\\[0.7ex]
$X \bigotimes Y$ is the Cartesian product of sets $X$ and $Y$.
\\[1.5ex]
\indent
It is not difficult to prove the following homomorphism theorem for $\Lambda $-factors in 3-connected graphs.
\bs 
\label{Gamma}
Let $B$ and each $A(u)$, 
$u \in V(B)$, be cubic 3-connected graphs.
Suppose that each $v(A(u))= 2 \bmod 6$ and 
each $(A(u), a^u)$ satisfies the following conditions: 
\\[1ex]
$(h1)$
$A(u) - (N^u \cup a^u \cup y)$ has no $\Lambda $-factor
for every vertex $y$ in $A(u) - (N^u \cup a^u)$ adjacent to a vertex in $N^u$
and
\\[1ex]
$(h2)$ $A(u) - \{a^u, z\}$ and $A(u) - W$
has a $\Lambda $-factor for every
$a^uz \in E(A(u))$ and a 5-vertex path $W$ in $A(u)$ centered at $a^u$, respectively.

Then
\\[1ex]
$(\gamma 1)$
$\Gamma (G, P) \cap \Gamma (G, Q) = \emptyset $ for 
$P, Q \in \Gamma (B)$, $P \ne Q$,
\\[1ex]
$(\gamma 2)$
$\Gamma (G, P) = 
( \bigotimes \{ \Gamma (A_1(v,P)): v \in V_1(P)\})
(\bigotimes \{ \Gamma (A_2(u,P)): u \in V_2(P)\}) $,
and
\\[1ex]
$(\gamma 3)$
$\Gamma (G) = \bigcup \{\Gamma (G, P): 
P \in \Gamma (B)\}$, 
\\[1ex] and so $G$ has a $\Lambda $-factor if and only if 
$B$ has a $\Lambda $-factor.
\es

By {\bf \ref{G,x,W}}, there are infinitely many pairs $(A,a)$  such that $A$ is  cubic 3-connected graph, 
$v(A) = 2 \bmod 6$,  $a \in V(A)$, and $(A,a)$ satisfies 
$(h1)$ in  {\bf \ref{Gamma}}.
If $(z1)$ is true, then by {\bf \ref{3-con}},
condition $(h2)$ in  {\bf \ref{Gamma}} is satisfied
for every pair $(A,a)$ such that $A$ is  cubic 3-connected graph, $v(A) = 2 \bmod 6$, and $a \in V(A)$.

\end{document}